\documentclass[final,hidelinks,onefignum,onetabnum]{siamart251216}

\usepackage{amssymb}
\usepackage{amsmath}
\usepackage{mathtools}
\usepackage{rotating}
\usepackage{latexsym, enumerate}
\usepackage{graphics}
\usepackage{epsfig}
\usepackage{tabularx}
\usepackage{eqnarray}
\usepackage{amsopn}
\usepackage{lipsum}

\usepackage{epstopdf}
\usepackage{algorithmic}
\ifpdf
  \DeclareGraphicsExtensions{.eps,.pdf,.png,.jpg}
\else
  \DeclareGraphicsExtensions{.eps}
\fi

\usepackage{fancybox}

\DeclareMathOperator{\sign}{sign}

\headers{A super-convergent method for advection-diffusion problems}{D. Y. Le Roux}

\title{A super-convergent quasi-discontinuous Galerkin method for advection-diffusion problems\thanks{Submitted to the editors DATE.
\funding{}}}

\author{Daniel Y. Le Roux\thanks{Universit\'e Lyon 1, CNRS, Institut Camille Jordan, 43, blvd du 11 novembre 1918, 69622 Villeurbanne Cedex, France
  (\email{dleroux@math.univ-lyon1.fr}).}}

\begin{document}

\maketitle

\begin{abstract}
The aim of the present study is to develop discontinuous Galerkin and nonconforming finite-element methods for solving the 2D advection-diffusion equation. The discrete formulations are stabilized by means of a family of numerical ﬂuxes and the dispersion relations are constructed via a Fourier/stability analysis. The linear discontinuous $P_{1}^{DG}$ and non conforming $P_{1}^{NC}$ Galerkin discretizations are studied. A long time stability result is proven for the $P_{1}^{NC}$ and $P_{1}^{DG}$ schemes and ﬂuxes examined here by inspecting the discrete frequencies, solutions of the dispersion relations. A super-convergent result is also demonstrated for the $P_{1}^{NC}$ discrete frequency ($O(h^4)$ accuracy) compared to the $P_{1}^{DG}$ one ($O(h^2)$ accuracy). These results are in line with those obtained in a previous paper for the shallow-water system. 
\end{abstract}



\section{Introduction}
Discontinuous Galerkin (DG) methods are a class of finite element methods using discontinuous basis functions, which are usually chosen as piecewise polynomials, to approximate the discrete solutions. The DG methods were first proposed by Reed and Hill \cite{Reed} in 1973 and LeSaint and Raviart \cite{LSR} in 1974, in the framework of neutron linear transport equation. Later, Cockburn {\it et al.} applied the DG method for solving nonlinear hyperbolic equations in a series of papers e.g. \cite{Cockburn6}. The DG methods have found their way into the main stream of computational fluid dynamics and other area of applications. For a detailed description of the method as well as its implementation and applications, we refer to the lecture notes \cite{cockburn, Shu} and the survey paper \cite{cockburn_survey}.

It is not straightforward to apply the DG method directly to equations containing higher order spatial derivatives, thereby the local DG (LDG) method was introduced. There are usually two steps to implement the LDG method: first rewriting the problem into a form which contains only first order derivatives, and secondly applying the DG method by choosing appropriate numerical fluxes to ensure stability. The first LDG method was constructed by Cockburn and Shu \cite{Cockburn_Shu} for solving time dependent convection diffusion systems, motivated by the work of Bassi and Rebay \cite{Bassy-Rebay} in the case of the compressible Navier-Stokes equations. Following these pioneering papers, LDG methods have been successfully designed and applied in a number of models involving diffusion and dispersive problems \cite{Castillo-2000, Cockburn}. The DG and LDG methods share several attractive properties, such as good stability, high order accuracy, and flexibility regarding $h-p$ adaptivity and on complex geometries. 

The aim of the present study is to develop high order discontinuous Galerkin and nonconforming finite-element methods for solving the 2D advection-diffusion equation. In \cite{LeRoux-JCP-2024}, a super-convergent result is demonstrated for the discrete frequencies solutions of the shallow-water system by using the non conforming $P_{1}^{NC}$ Galerkin discretization. This is why the linear non conforming $P_{1}^{NC}$ and discontinuous $P_{1}^{DG}$ Galerkin discretizations are studied and compared in this paper. The continuous model is presented in Section~2 and the discretization of the model equations is performed in Section~3 along with the introduction of a family of numerical ﬂuxes. The Fourier/dispersion analyses are then conducted in Section~4, and the results are analysed. A long time stability result is proven for the $P_{1}^{NC}$ and $P_{1}^{DG}$ schemes and ﬂuxes examined here by inspecting the discrete frequencies, solutions of the dispersion relations. A super-convergent result is also demonstrated for the $P_{1}^{NC}$ discrete frequency ($O(h^4)$ accuracy) compared to the $P_{1}^{DG}$ one ($O(h^2)$ accuracy). Some concluding remarks complete the study in Section~5.
\section{The problem model}
\label{Problem-model}
Let $\Omega$ be a Lipschitz domain in $ \mathbb{R}_{}^{2}$ with boundary $\partial \Omega = \Gamma_{D}^{} \cup \Gamma_{N}^{}$ satisfying $\Gamma_{D}^{} \cap \Gamma_{N}^{} = \emptyset$, where $\Gamma_{D}^{}$ and $\Gamma_{N}^{}$ are introduced in order to be able to handle Dirichlet ($D$) and Neumann ($N$) boundary conditions. For a given finite time $T>0$, consider the advection-diffusion problem
\begin{alignat}{2}
\label{cont-eq}
\frac{\partial u}{\partial t} + \nabla \cdot (u\, \mathbf{a}) - \nu \Delta u & = f, \qquad && \mbox{ in } \Omega \times (0,T), \\
\label{bnd_D}
u & = u_{D}^{}, \qquad && \mbox{ on } \Gamma_{D}^{} \times (0,T), \\
\label{bnd_N}
\frac{\partial u}{\partial n} & = \boldsymbol{\theta}_{N}^{} \cdot \mathbf{n}, \qquad && \mbox{ on } \Gamma_{N}^{} \times (0,T), \\
\label{cont-init}
u(\mathbf{x},t=0) & = u_{0}^{}(\mathbf{x}), \qquad && \mbox{ in } \Omega,
\end{alignat}
where $\mathbf{a}=(a_{1}^{},a_{2}^{})$ is the constant advection speed of the solution with $\mid \! \mathbf{a} \! \mid = (a_{1}^{2}+a_{2}^{2})_{}^{\frac{1}{2}}$, $\nu$ is the constant diffusion coefficient, $\mathbf{x}=(x,y)$ and $f=f(\mathbf{x},t)$ is a source term which introduces energy into the system. Further, we consider $a_{1}^{}>0$ and $a_{2}^{}>0$ without loss of generality. The initial solution $u_{0}^{}(\mathbf{x})$ and $f(\mathbf{x},t)$ are chosen in $L_{}^{2}(\Omega)$, $\mathbf{n}$ is the unit outward normal vector to $\partial \Omega$ and we assume that the given functions $u_{D}^{}$ and $\boldsymbol{\theta}_{N}^{}$ are sufficiently regular on $\Gamma_{D}^{} \times (0,T)$ and $\Gamma_{N}^{} \times (0,T)$, respectively. 


To establish the \textit{local} $P_{1}^{NC}$ and $P_{1}^{DG}$ methods, the flux $\boldsymbol{\theta}$ is introduced as an auxiliary variable, and~\crefrange{cont-eq}{cont-init} are rewritten as a system of first-order equations
\begin{alignat}{2}
\label{cont-eq-2}
\boldsymbol{\theta} & = \nu \nabla u, \qquad && \mbox{ in } \Omega \times (0,T), \\
\label{cont-eq-1}
\frac{\partial u}{\partial t} + \nabla \cdot (u\, \mathbf{a}) - \nabla \cdot \boldsymbol{\theta} & = f,  \qquad && \mbox{ in } \Omega \times (0,T), \\
\label{bnd_D_2}
u & = u_{D}^{}, \qquad && \mbox{ on } \Gamma_{D}^{} \times (0,T), \\
\label{bnd_N_2}
\boldsymbol{\theta} \cdot \mathbf{n} & = \boldsymbol{\theta}_{N}^{} \cdot \mathbf{n}, \qquad && \mbox{ on } \Gamma_{N}^{} \times (0,T), \\
\nonumber u(\mathbf{x},t=0) & =  u_{0}^{}(\mathbf{x}). \qquad && \mbox{ in } \Omega.
\end{alignat}
In order to describe the weak formulation, a few notations are first defined. Let $\{\tau_{h}^{}\}^{}_{h>0}$ denote a partition of the domain $\Omega$ into a finite number $N^{}_{el}$ of disjoint open non degenerate elements $K_{el}^{}$, $el=1,2,\cdots, N^{}_{el}$. We denote by $h_{K_{el}^{}}^{}$ the diameter of $K_{el}^{}$ and we set $h := \max_{K_{el}^{} \in \tau_{h}^{}}^{} K_{el}^{}$. The closure of $\Omega$ satisfies $\overline{\Omega} = \cup^{N^{}_{el}}_{el=1} \overline{K}_{el}^{}$ and for $K_{el_{}^{+}}, K_{el_{}^{-}} \in \tau_{h}^{}$ we have $\overset{\circ}{K}_{el_{}^{+}} \cap \overset{\circ}{K}^{}_{el_{}^{-}} = \emptyset$ if $ K_{el_{}^{+}} \neq K_{el_{}^{-}}$, where $\overset{\circ}{K}_{el}^{}$ is the interior of $K_{el}^{}$. In the case $K_{el_{}^{+}} \cap K^{}_{el_{}^{-}} \neq \emptyset$, then  $K_{el_{}^{+}} \cap K^{}_{el_{}^{-}}$ is either a common face or a common vertex of $K_{el_{}^{+}}$ and $K_{el_{}^{-}}$. 

Each element $K_{el}^{}$ is supposed to have a Lipschitz boundary $\partial K_{el}^{}$. Let $\Gamma$ be the finite ensemble of $M^{}_{ed}$ inter-element boundaries denoted by $\Gamma_{ed}^{} = \overline{\partial K^{}}_{el_{}^{+}} \cap \overline{\partial K^{}}_{el_{}^{-}}$, $ed=1,2,\cdots, M^{}_{ed}$, where $\Gamma_{ed}^{}$ is an interior edge (not belonging to $\partial \Omega$) shared by elements $K_{el_{}^{+}}$ and $K^{}_{el_{}^{-}}$, such that
\[\overline{\Gamma} = \cup^{M^{}_{ed}}_{ed=1} \overline{\Gamma}^{}_{ed} \quad \mbox{and} \quad
\Gamma^{}_{ed_{}^{+}} \cap \Gamma^{}_{ed_{}^{-}} = \emptyset \mbox{\; if \;} ed_{}^{+} 
\neq ed_{}^{-}.\]
Each inter-element boundary $\Gamma^{}_{ed} \in \Gamma$ is associated with a unit normal vector denoted by $\mathbf{n}_{ed}^{}=({n}^{x}_{ed},{n}^{y}_{ed})$ at $\mathbf{x} \in \Gamma^{}_{ed}$, and $\mathbf{n}_{ed}^{+}$ and $\mathbf{n}_{ed}^{-}$ are the unit normal vectors on $\Gamma_{ed}^{}$ pointing exterior to elements $K_{el_{}^{+}}$ and $K^{}_{el_{}^{-}}$, respectively. 
\section{Discretization of the model equations}
\label{numerical:heat}
\subsection{The discrete weak formulation}
Let $\mathcal{P}_{}^{r} (K_{el}^{})$ be the space of polynomials of degree at most $r \geq 0$, defined on $K_{el}^{}$, and consider the polynomial spaces $U_{h,K_{el}^{}}^{}$ and $\boldsymbol{\Theta}^{}_{h,K_{el}^{}}$ of piecewise smooth functions that are differentiable over an element, but which allows discontinuities between elements. These spaces will be defined in~\cref{NC-discrete,DG-discrete} for the $P_{1}^{NC}$ and $P_{1}^{DG}$ schemes, respectively. The components of $u$ and $\boldsymbol{\theta}$ are approximated in $U_{h,K_{el}^{}}^{}$ and $\boldsymbol{\Theta}^{}_{h,K_{el}^{}}$, respectively, and these are designated by $u_{h}^{}$ and $\boldsymbol{\theta}_{h}^{}$. 

The weak formulation is obtained by multiplying~\eqref{cont-eq-1} and~\eqref{cont-eq-2} by arbitrary test functions $v(\mathbf{x}) \in U^{}_{h,K_{el}^{}}$ and $\boldsymbol{\psi} (\mathbf{x}) \in \boldsymbol{\Theta}^{}_{h,K_{el}^{}}$, respectively, and the x- or y-component of $\boldsymbol{\psi} (\mathbf{x})$ is formally denoted by $\psi (\mathbf{x})$. Integrating the flux terms by parts over element $K_{el}^{}$, using Green's theorem, leads to find $u \in U^{}_{h,K_{el}^{}}$ and $\boldsymbol{\theta} \in \boldsymbol{\Theta}^{}_{h,K_{el}^{}}$ satisfying
\begin{align}
\label{WF-eq-1}
& \int^{}_{K_{el}^{}} \frac{\partial u_{h}^{}}{\partial t} \, v \, d\mathbf{x} = \int^{}_{K_{el}^{}} u_{h}^{}\, \mathbf{a} \cdot \nabla v \, d\mathbf{x} - \int^{}_{\partial K_{el}^{}} u_{h}^{}\, \mathbf{a} \cdot \mathbf{n} ^{}_{ed}\, v \, ds 
+ \int^{}_{K_{el}^{}} f \, v \, d\mathbf{x} \\ 
& \hspace{2.4cm} - \int^{}_{K_{el}^{}} \boldsymbol{\theta}_{h}^{} \cdot 
\nabla v \, d\mathbf{x} + \int^{}_{\partial K_{el}^{}}  \boldsymbol{\theta}_{h}^{} \cdot \mathbf{n} ^{}_{ed}\, v \, ds, \nonumber \hspace{0.9cm} \forall v \in U^{}_{h,K_{el}^{}}, \\
\label{WF-eq-2}
& \int^{}_{K_{el}^{}} \boldsymbol{\theta}_{h}^{} \cdot \boldsymbol{\psi} \, d\mathbf{x} =  - \nu \int^{}_{K_{el}^{}} u_{h}^{} \, \nabla \cdot \boldsymbol{\psi} \, d\mathbf{x} + \nu \int^{}_{\partial K_{el}^{}} u_{h}^{} \, \boldsymbol{\psi} \cdot \mathbf{n}^{}_{ed} \, ds, \quad \, \forall \boldsymbol{\psi} \in \boldsymbol{\Theta}^{}_{h,K_{el}^{}}.
\end{align}
%
%
%
%
%
%
Then we sum~\eqref{WF-eq-1} and~\eqref{WF-eq-2} over all elements $K_{el}^{} \in \tau_{h}^{}$. As $u_{h}^{}$ and $\boldsymbol{\theta}_{h}^{}$ may be discontinuous across inter-element edges $\Gamma^{}_{ed} \in \Gamma^{}$, they must be replaced by $(u_{h}\, \mathbf{a})_{}^{\star}$, $u_{h}^{*}$ and $\boldsymbol{\theta}_{h}^{*}$ which denote the numerical traces of $u_{h}^{}\, \mathbf{a}$, $u_{h}^{}$ and $\boldsymbol{\theta}_{h}^{}$ on $\Gamma^{}_{ed}$. The approximate solution ($u_{h}^{}$, $\boldsymbol{\theta}_{h}^{}$) belonging to $U^{}_{h,K_{el}^{}} \times \boldsymbol{\Theta}^{}_{h,K_{el}^{}}$ is then obtained from~\eqref{WF-eq-1} and~\eqref{WF-eq-2} as the solution of
\begin{align}
\label{WFh-eq-1-0}
 \sum_{el=1}^{N^{}_{el}} \int^{}_{K_{el}^{}} \frac{\partial u_{h}^{}}{\partial t} \, v \, d\mathbf{x} & = \sum_{el=1}^{N^{}_{el}} \int^{}_{K_{el}^{}} u_{h}^{}\, \mathbf{a} \cdot \nabla v \, d\mathbf{x} - \sum_{el=1}^{N^{}_{el}} \int^{}_{\partial K_{el}^{}} (u_{h}\, \mathbf{a})_{}^{\star} \cdot\mathbf{n}^{}_{ed}\,  v \, ds \\ 
 + \sum_{el=1}^{N^{}_{el}} \int^{}_{K_{el}^{}} f \, v \, d\mathbf{x} - \sum_{el=1}^{N^{}_{el}} \int^{}_{K_{el}^{}} & \boldsymbol{\theta}_{h}^{} \cdot 
\nabla v \, d\mathbf{x} + \sum_{el=1}^{N^{}_{el}} \int^{}_{\partial K_{el}^{}} \boldsymbol{\theta}_{h}^{*} \cdot \mathbf{n} ^{}_{ed}\, v \, ds, \quad \forall v \in U^{}_{h,K_{el}^{}}, \nonumber \\
\label{WFh-eq-2-0}
 \sum_{el=1}^{N^{}_{el}} \int^{}_{K_{el}^{}} \boldsymbol{\theta}_{h}^{} \cdot \boldsymbol{\psi} \, d\mathbf{x} & = -\nu \sum_{el=1}^{N^{}_{el}} \int^{}_{K_{el}^{}} u_{h}^{} \, \nabla \cdot \boldsymbol{\psi} \, d\mathbf{x}  \\ 
& + \nu \sum_{el=1}^{N^{}_{el}} \int^{}_{\partial K_{el}^{}} u_{h}^{*} \, \boldsymbol{\psi} \cdot \mathbf{n}^{}_{ed} \, ds, \quad \forall \boldsymbol{\psi} \in \boldsymbol{\theta}^{}_{h,K_{el}^{}}, \nonumber
\end{align}
To complete the definition of the DG method in~\eqref{WFh-eq-1-0} and~\eqref{WFh-eq-2-0}, it remains to choose $(u_{h}\, \mathbf{a})_{}^{\star}$, $\boldsymbol{\theta}_{h}^{*}$ and $u_{h}^{*}$ on inter-element boundaries $\Gamma_{ed}^{}$ as to render the method consistent and stable. The numerical trace is a \textit{single valued} function defined on $\Gamma^{}_{ed}$ and it usually depends on the numerical solution from both sides of the interface using the mean values and jumps at $\mathbf{x} \in \Gamma_{ed}^{}$. For any function $\chi \in U^{}_{K_{el}^{}}$ or $\boldsymbol{\theta}^{}_{K_{el}^{}}$, where $\chi$ stands for $u$, $\theta_{}^{x}$ or $\theta_{}^{y}$, with $\boldsymbol{\theta}=(\theta_{}^{x},\theta_{}^{y})$, and for each element $K_{el}^{}$, the left and right values of $\chi$ on $\Gamma^{}_{ed}$ from the interior of $K_{el_{}^{\pm}}^{}$ are represented by $\chi^{\pm}_{}$ with
\begin{equation}
\chi^{-}_{}(\mathbf{x}) = \lim^{}_{\epsilon \to 0^{-}_{}} 
\chi(\mathbf{x} +  \epsilon \mathbf{n}_{ed}^{}), \quad 
\chi^{+}_{}(\mathbf{x}) = \lim^{}_{\epsilon \to 0^{+}_{}} \chi(\mathbf{x} + 
\epsilon \mathbf{n}_{ed}^{}). 
\label{defiLiR}
\end{equation}
The mean values $\{\, \cdot \, \}$ and jumps $[\![\,  \cdot \,]\!]$ at $\mathbf{x} \in \Gamma^{}_{ed}$ are then defined as 
\begin{align*}
 \{ u \}  &:= \tfrac{1}{2} \left(u^{+}_{} + u^{-}_{} \right),  &
 \{ \boldsymbol{\theta} \}  &:= \tfrac{1}{2} \left(\boldsymbol{\theta}^{+}_{} + \boldsymbol{\theta}^{-}_{} \right), \\
 [\![ u ]\!] \,&:= u^{+}_{} \,\mathbf{n}_{ed}^{+} + u^{-}_{} \,\mathbf{n}_{ed}^{-}, &
 [\![ \boldsymbol{\theta} ]\!] \,&:= \boldsymbol{\theta}^{+}_{} \cdot \mathbf{n}_{ed}^{+} + \boldsymbol{\theta}^{-}_{} \cdot \mathbf{n}_{ed}^{-}.
\end{align*}
The advantage of the above definitions is that they do not depend on assigning an ordering to the elements $K_{el}^{}$.

In this study the so-called classical up-winding numerical flux is chosen for $(u_{h}\, \mathbf{a})_{}^{\star}$ while, following \cite{Castillo-2000},  general numerical traces are introduced for $u_{h}^{\star}$ and $\boldsymbol{\theta}_{h}^{*}$ as follows
\begin{align}
(u_{h}\, \mathbf{a})_{}^{\star} & = \{u_{h}^{} \}\,\mathbf{a} + \frac{1}{2} \mid \! \mathbf{a} \cdot \mathbf{n}_{ed}^{} \! \mid \, [\![ u_{h}^{} ]\!], \label{tr_uhc}\\
\boldsymbol{\theta}_{h}^{*} & = \{ \boldsymbol{\theta}_{h}^{} \} - \boldsymbol{\gamma} \, [\![ \boldsymbol{\theta}_{h}^{} ]\!] 
- \zeta h_{K_{el}^{}}^{m_{1}} \,[\![ u_{h}^{} ]\!], \label{tr_qh} \\  
u_{h}^{*} & = \{ u_{h}^{} \} + \boldsymbol{\gamma} \cdot [\![ u_{h}^{} ]\!]
- \xi h_{K_{el}^{}}^{m_{2}} \, [\![ \boldsymbol{\theta}_{h}^{} ]\!], \label{tr_uh}
\end{align}
with $\boldsymbol{\gamma} \cdot \mathbf{n}^{\pm}_{ed} = \frac{1}{2} \sign (\mathbf{b} \cdot \mathbf{n}^{\pm}_{ed})$, where $\mathbf{b}$ is the arbitrary chosen vector $\mathbf{b} = \frac{1}{\sqrt{2}} (1,1)_{}^{t}$. The role of $\mathbf{b}$ is to provide a single rule to pick the numerical traces $u_{h}^{*}$ and $\boldsymbol{\theta}_{h}^{*}$. 

Note that $\zeta$ and $\xi$ are assumed to be constant parameters, while the coefficients $h_{K_{el}^{}}^{m_{1}}$ and $h_{K_{el}^{}}^{m_{2}}$ being mesh-dependent, allow the possibility that a suitable choice of $m_{1}$ and $m_{2}$ could improve the order of convergence of the method. However, in order to satisfy homogeneity constraints in~\eqref{tr_qh} -~\eqref{tr_uh}, we impose $m_{1}=-1$ and $m_{2}=~1$. Further, after obtaining the dispersion relations in~\cref{Disp-Rel}, it was observed that $\xi$ does not provide any particular information, both in terms of accuracy and stability and we let $\xi=0$. In this case, the auxiliary variable $\boldsymbol{\theta}_{h}^{}$ can be locally solved in terms of $u_{h}^{}$ by using~\eqref{WFh-eq-2-0} and then easily eliminated from~\eqref{WFh-eq-1-0}. This \textit{local} and computationally efficient property gives its name to the LDG method. In addition, we let $\zeta h_{K_{el}^{}}^{-1} = \dfrac{p}{h^{}_{ed}}$ in the following, where $p$ is a user-dependent real parameter and $h^{}_{ed}$ is equal to the diameter of the inter-element boundary $\Gamma^{}_{ed}$, i.e. a local length scale, which yields
\begin{alignat}{2}
\label{jump_q}
\boldsymbol{\theta}_{h}^{*} & = \{\boldsymbol{\theta}_{h}^{} \} - \boldsymbol{\gamma} \, [\![ \boldsymbol{\theta}_{h}^{} ]\!] 
- \frac{p}{h^{}_{ed}} \,[\![ u_{h}^{} ]\!], \quad & \mbox{ on } \Gamma^{}_{ed}, \\  
\label{jump_u}
u_{h}^{*} & = \{ u_{h}^{} \} + \boldsymbol{\gamma} \cdot [\![ u_{h}^{} ]\!], \quad & \mbox{ on } \Gamma^{}_{ed}.
\end{alignat}
and following \cite{Cockburn, Riviere} we have
\begin{alignat}{2}
\label{jump_q_2}
\boldsymbol{\theta}_{h}^{*} & = \left\{ 
\begin{array}{l}
\boldsymbol{\theta}_{h}^{+}-\dfrac{p}{h^{}_{ed}} \,(u_{h}^{+}-u_{D}^{}) \mathbf{n}_{ed}^{+}, \quad  \mbox{ on } \Gamma^{}_{D}, \\*[0.2cm]
\boldsymbol{\theta}_{N}^{}, \hspace*{3.15cm} \quad  \mbox{ on } \Gamma^{}_{N},
\end{array} \right.
\\ 
\label{jump_u_2}
u_{h}^{*} & = 
\left\{ 
\begin{array}{l}
u_{D}^{}, \quad  \mbox{ on } \Gamma^{}_{D}, \\*[0.2cm]
u_{h}^{+}, \quad  \mbox{ on } \Gamma^{}_{N}.
\end{array} \right.
\end{alignat}
Interior edges $\Gamma^{}_{ed}$ appear twice in the sum over element boundary integrals on $\partial K_{el}^{}$ in~\eqref{WFh-eq-1} and~\eqref{WFh-eq-2}. Line integral contributions hence appear from both side of the neighbouring elements with the normal $ \mathbf{n}_{ed}^{}$ pointing in the opposite direction in each element, revealing the jumps of the test functions. By using the numerical traces~\eqref{jump_q}~-~\eqref{jump_u_2}, for all $v \in U^{}_{h,K_{el}^{}}$ and for all $\boldsymbol{\psi} \in \boldsymbol{\theta}^{}_{h,K_{el}^{}}$ in~\eqref{WFh-eq-1} and~\eqref{WFh-eq-2}, repectively, we obtain
\begin{align}
\label{WFh-eq-1}
 \sum_{el=1}^{N^{}_{el}} \int^{}_{K_{el}^{}} \frac{\partial u_{h}^{}}{\partial t} \, v \, d\mathbf{x} 
 & = \sum_{el=1}^{N^{}_{el}} \int^{}_{K_{el}^{}} u_{h}^{}\, \mathbf{a} \cdot \nabla v \, d\mathbf{x} - \sum_{el=1}^{N^{}_{el}} \int^{}_{\partial K_{el}^{}} (u_{h}\, \mathbf{a})_{}^{\star} \cdot\mathbf{n}^{}_{ed}\,  v \, ds \\ 
+ \sum_{el=1}^{N^{}_{el}} \int^{}_{K_{el}^{}} f \, v \, d\mathbf{x} - & \sum_{el=1}^{N^{}_{el}} \int^{}_{K_{el}^{}} 
 \boldsymbol{\theta}_{h}^{} \cdot \nabla v \, d\mathbf{x} \nonumber  \\
+ \nu \sum_{ed=1}^{M^{}_{ed}} \int^{}_{\Gamma_{ed}^{}} \Big( \{\boldsymbol{\theta}_{h}^{} \} & - \boldsymbol{\gamma} \, [\![ \boldsymbol{\theta}_{h}^{} ]\!] 
- \frac{p}{h^{}_{ed}} \,[\![ u_{h}^{} ]\!] \Big) \cdot [\![ v ]\!]\,ds \nonumber  \\
+ \nu \int^{}_{\Gamma_{D}^{}} \Big( \boldsymbol{\theta}_{h}^{+} - \dfrac{p}{h^{}_{ed}} & \,(u_{h}^{+} - u_{D}^{}) \mathbf{n}_{ed}^{+} \Big) \cdot \mathbf{n}^{+}_{ed} \,v\, ds 
+ \nu \int^{}_{\Gamma_{N}^{}} \boldsymbol{\theta}_{N}^{} \cdot \mathbf{n}^{+}_{ed} \,v\, ds, 
 \nonumber \\
\label{WFh-eq-2}
 \sum_{el=1}^{N^{}_{el}} \int^{}_{K_{el}^{}} \boldsymbol{\theta}_{h}^{} \cdot \boldsymbol{\psi} \, d\mathbf{x} & = -\nu \sum_{el=1}^{N^{}_{el}} \int^{}_{K_{el}^{}} u_{h}^{} \, \nabla \cdot \boldsymbol{\psi} \, d\mathbf{x} \\
+ \nu \sum_{ed=1}^{M^{}_{ed}} \int^{}_{\Gamma_{ed}^{}} \Big( \{ u_{h}^{} \} + \boldsymbol{\gamma} \cdot &[\![ u_{h}^{} ]\!] \Big) [\![ \boldsymbol{\psi}_{h}^{} ]\!]\,ds 
+  \nu \int^{}_{\Gamma_{D}^{}} u_{D}^{}  \boldsymbol{\psi} \cdot \mathbf{n}^{+}_{ed} \, ds 
+  \nu \int^{}_{\Gamma_{N}^{}} u_{h}^{+} \boldsymbol{\psi} \cdot \mathbf{n}^{+}_{ed} \, ds. \nonumber
\end{align}
The discretized weak formulation~\eqref{WFh-eq-1} for $u_{h}^{}$ and~\eqref{WFh-eq-2} for $\boldsymbol{\theta}_{h}^{}$ are computed in~\cref{NC-discrete,DG-discrete} using the $P^{NC}_{1}$ and $P^{DG}_{1}$ schemes, respectively. For the purpose of the subsequent Fourier analyses in~\cref{Four-disp}, a uniform mesh made up of biased right triangles is considered, denoted as Mesh~1, and the representative mesh length parameter $h$ is thus taken as a constant in the $x$- and $y$-directions. Further, periodic boundary conditions are considered and we let $f(\mathbf{x},t)=0$. In addition, only a few set of discrete equations have to be computed at typical nodal locations due to symmetry reasons and the needs of the Fourier analysis. 

Equations~\cref{WFh-eq-1,WFh-eq-2} are computed over triangle $K^{}_{1}$ (right panel) in~\cref{Node-distrib-P1NC} at nodes $1,2,3,$ for the $P^{NC}_{1}$ scheme, and over triangles $K^{}_{1}$ and $K^{}_{2}$ (right panel) in~\cref{Node-distrib-P1DG} at nodes $1,2,3,4,5,6,$ for the $P^{DG}_{1}$ approach. At such \textit{selected} nodes, and for both schemes, equations~\cref{WFh-eq-1,WFh-eq-2} may be expressed on the compact form
\begin{align}
\label{compact-discrete-eqn-u}
\frac{\partial }{\partial t}\left(u^{}_{1}, u^{}_{2},..., u^{}_{n}\right)_{}^{t} \! & = \!
-\frac{a_{1}^{}}{h} \boldsymbol{\Lambda}^{}_{1}(\{u^{}_{j}\}) - \frac{a_{2}^{}}{h} \boldsymbol{\Lambda}^{}_{2}(\{u^{}_{j}\})\! +\! \frac{p}{h^{2}_{}} \boldsymbol{\Lambda}^{}_{3}(\{u^{}_{j}\})\! + \!\frac{1}{h} \boldsymbol{\Lambda}^{}_{4}(\{\boldsymbol{\theta}^{}_{j}\}), \\*[0.2cm]
\label{compact-discrete-eqn-qx}
(\theta^{x}_{1}, \theta^{x}_{2}, ..., \theta^{x}_{n})_{}^{\,t} & = \frac{\nu}{h} \boldsymbol{\Lambda}^{}_{1}(\{u^{}_{j}\}), \\*[0.2cm]
\label{compact-discrete-eqn-qy}
(\theta^{y}_{1}, \theta^{y}_{2}, ..., \theta^{y}_{n})_{}^{\,t} & = \frac{\nu}{h} \boldsymbol{\Lambda}^{}_{2}(\{u^{}_{j}\}),
\end{align}
where $j=1,2,3,..,$ denotes a mesh node number. A total of $3n$ selected discrete equations are hence computed for the $P^{NC}_{1}$ (with $n=3$) and $P^{DG}_{1}$  (with $n=6$) schemes. The vectors $\boldsymbol{\Lambda}^{}_{1}, \boldsymbol{\Lambda}^{}_{2}, \boldsymbol{\Lambda}^{}_{3}, \boldsymbol{\Lambda}^{}_{4},$ are computed in~\cref{NC-discrete} at node $2$ of~\cref{Node-distrib-P1NC} for the $P^{NC}_{1}$ scheme and in~\cref{DG-discrete} at node $2$ of~\cref{Node-distrib-P1DG} for the $P^{DG}_{1}$ method, as examples. However, general results at the \textit{selected} nodes are given in~\cref{app-def-submatNC,app-def-submatDG} for the $P^{NC}_{1}$ and $P^{DG}_{1}$ discretizations, respectively.
\subsection{The \boldmath$P_{1}^{NC}$ discrete equations}
\label{NC-discrete}
\begin{figure}[h!]
\begin{center}
\begin{tabular}{lll}
\includegraphics*[height=34mm]{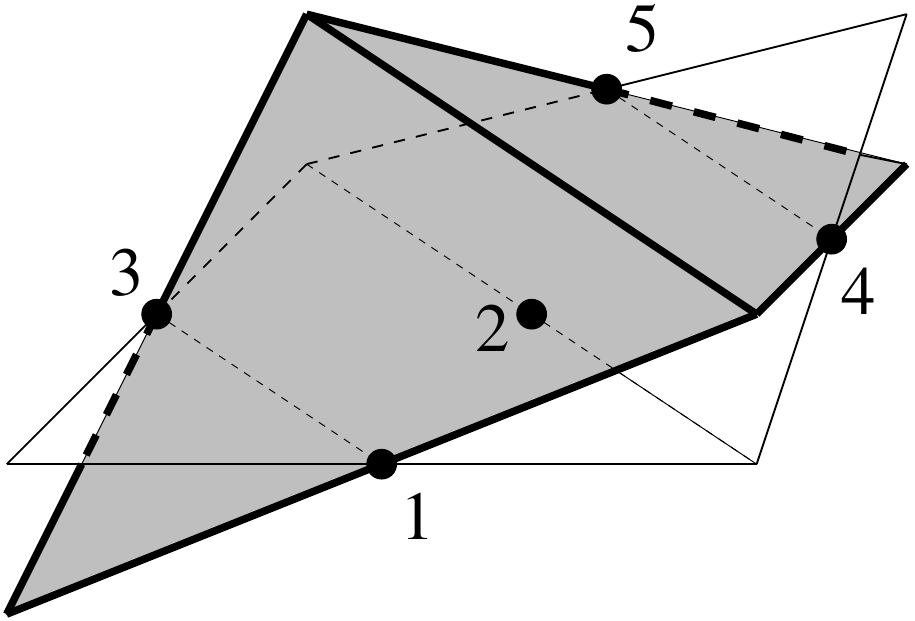}
& &
\includegraphics*[height=42mm]{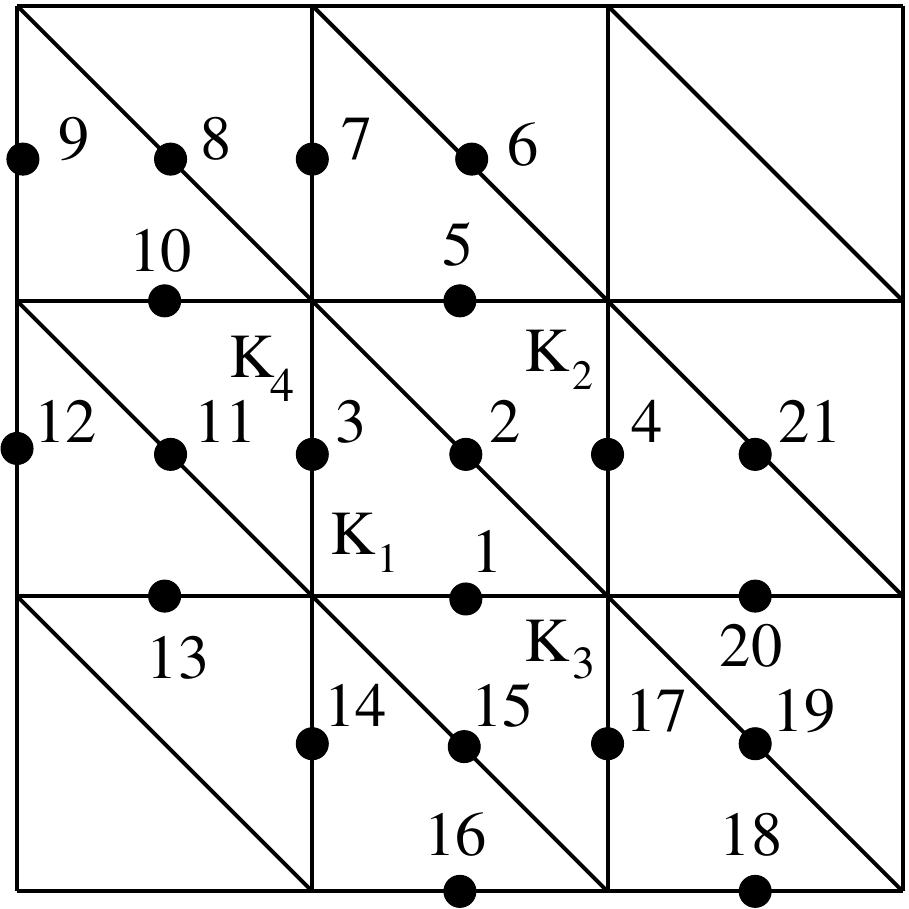}
\end{tabular}
\caption{The $P_{1}^{NC}$ linear basis function at node $2$ in grey (left). Node locations, used to compute the discrete equations, are represented by the symbol {\Large $\bullet$} (right).}
\label{Node-distrib-P1NC}
\end{center}
\end{figure}

The discrete solutions $u_{h}^{}$ and $\boldsymbol{\theta}_{h}^{}$, as well as the basis functions, belong to the space of polynomials $\mathcal{P}_{}^{1} (K_{el}^{})$, and hence they are approximated linearly over triangles $K_{el}^{}$ of $\tau_{h}^{}$. However, their nodal values are located at midpoints of edges of $K_{el}^{}$, as shown in~\cref{Node-distrib-P1NC} (right panel), and two adjacent triangles share the \textit{same} nodal value at the midpoint of their common edge. 

The shape of the $P_{1}^{NC}$ linear basis function that is used to approximate the components of $u_{h}^{}$ and $\boldsymbol{\theta}_{h}^{}$ on the element's two-triangle support is represented in grey at node $2$ in~\cref{Node-distrib-P1NC} (left panel). The basis function has the value $1$ at the midpoint and along the common edge of the two adjacent triangles and $0$ at the midpoints of the remaining four edges, represented by the symbol {\Large $\bullet$} in~\cref{Node-distrib-P1NC}. Since such a representation is only continuous across triangle boundaries at midpoint of the edges, and discontinuous everywhere else around a triangle boundary, this element is termed nonconforming ($NC$) in the FE literature. 

Consequently, the polynomial space $V_{h,K_{el}^{}}^{}$ is defined as
\begin{align*}
V_{h,K_{el}^{}}^{} := V_{h,K_{el}^{}}^{NC} = & \{u_{h}^{} : \, u_{h}^{} \vert^{}_{K_{el}^{}} \in 
\mathcal{P}_{}^{1} (K_{el}^{}) \mbox{ for all } K_{el}^{} \in \tau_{h}^{},  \\*[-0.1cm]
 & u_{h}^{} \mbox{ is continuous at the midpoints of $\Gamma_{ed}^{}$} \}.
\end{align*}
A useful property of the $P^{NC}_{1}$ element is that the mass matrix is diagonal since the $P_{1}^{NC}$ basis functions are orthogonal \cite{Thomasset-1981}, a desirable property that enhances computational efficiency.

The equation for $u_{h}^{}$ in~\eqref{WFh-eq-1} is computed at node $2$ of~\cref{Node-distrib-P1NC} (right panel), and the compact support of the NC basis function $v^{}_{2}$ at node $2$ is now made up of triangles $K^{}_{1}$ and $K^{}_{2}$, namely $K^{}_{1} \cup K^{}_{2}$, leading to
\begin{align}
\label{WFh-eq-11}
\int^{}_{K^{}_{1} \cup K^{}_{2}} \frac{\partial u_{h}^{}}{\partial t} \, v^{}_{2} \, d\mathbf{x} & = \int^{}_{K^{}_{1} \cup K^{}_{2}} u_{h}^{}\, \mathbf{a} \cdot \nabla v^{}_{2} \, d\mathbf{x} - \int^{}_{\partial K^{}_{1} \cup K^{}_{2}} (u_{h}\, \mathbf{a})_{}^{\star} \cdot \mathbf{n}^{}_{ed}\,  v^{}_{2} \, ds \\ 
& - \int^{}_{K^{}_{1} \cup K^{}_{2}} \boldsymbol{\theta}_{h}^{} \cdot 
\nabla v^{}_{2} \, d\mathbf{x} + \int^{}_{\partial K^{}_{1} \cup K^{}_{2}} \boldsymbol{\theta}_{h}^{*} \cdot \mathbf{n}^{}_{ed}\, v^{}_{2} \, ds. \nonumber 
\end{align}
The integrals in~\eqref{WFh-eq-11} are computed by expanding $u_{h}^{}$ and $\boldsymbol{\theta}_{h}^{}$ over triangles $K_{1}^{}$ and $K_{2}^{}$ in terms of the linear NC basis functions $v_{j}^{}$ and $\psi^{}_{j}, j=1,2,3,4,5$, respectively. For example, $u_{h}^{} \vert^{}_{K_{1}^{}}= u^{}_{1} v^{}_{1} + u^{}_{2} v^{}_{2} + u^{}_{3} v^{}_{3}$. Let $\varphi$ represent the classical conforming linear basis over $K_{el}$ which takes the value $1$ at a given vertex of $K_{el}^{}$ and the value $0$ at the remaining two vertices of $K_{el}^{}$. The basis functions $\varphi^{}_{r}$ and $v^{}_{q}$ are linked over $K^{}_{el}$ by the relation $v^{}_{q}=1-2\varphi^{}_{r}$ if $q$ is the mid-edge point of $K^{}_{el}$ facing the vertex $r$ of $K^{}_{el}$. Due to the orthogonality property of the NC basis functions, with $\int_{\Omega} v^{}_{i} v^{}_{j}\, d\Omega=\frac{S^{}_{j}}{3}\delta^{}_{ij}$, where $S^{}_{j}$ is the area of the support of $v^{}_{j}$ and $\delta^{}_{ij}$ is the Kronecker delta, and elementary computations leading to $\frac{\partial v^{}_{2}}{\partial x}|^{}_{K_1} = \frac{\partial v^{}_{2}}{\partial y}|^{}_{K_1} = -\frac{\partial v^{}_{2}}{\partial x}|^{}_{K_2} = -\frac{\partial v^{}_{2}}{\partial y}|^{}_{K_2}= \frac{2}{h}$, we obtain
\begin{align*}
\int^{}_{K^{}_{1} \cup K^{}_{2}} \frac{\partial u_{h}^{}}{\partial t}\,v^{}_{2}\,d{\mathbf{x}}
& = \frac{h_{}^{2}}{3} \frac{\partial u_{2}^{}}{\partial t}, \\
\int^{}_{K^{}_{1} \cup K^{}_{2}} u_{h}^{}\, \mathbf{a} \cdot \nabla v^{}_{2} \, d\mathbf{x} 
& = \frac{h}{3} (a_{1}^{}+a_{2}^{}) (u_{1}^{}+u_{3}^{}-u_{4}^{}-u_{5}^{}) \\
- \int^{}_{K^{}_{1} \cup K^{}_{2}} \boldsymbol{\theta}_{h}^{} \cdot \nabla v^{}_{2} \, d\mathbf{x} 
& = \frac{h}{3} \big(-\theta^{x}_{1} -\theta^{x}_{3} + \theta^{x}_{4} + \theta^{x}_{5} - \theta^{y}_{1} -\theta^{y}_{3} +\theta^{y}_{4} +\theta^{y}_{5} \big).
\end{align*}
Let $\Gamma_{j}^{}$ and $\mathbf{n}_{j}^{}$ be the edge corresponding to midpoint $j$ and the associated normal, respectively, and let $u^{}_{\,\Gamma_{i}^{}}|^{}_{K^{}_{j}}$ denote the evaluation of $u^{}$ along $\Gamma_{i}^{}$ on $K^{}_{j}$. Consequently, $\partial K^{}_{1} = \Gamma_{1}^{} \cup \Gamma_{2}^{}\cup \Gamma_{3}^{}$ and $\partial K^{}_{2} = \Gamma_{2}^{} \cup \Gamma_{4}^{}\cup \Gamma_{5}^{}$. 
%

Along $\Gamma_{1}^{}$, for example, we have $\mathbf{n}_{1}^{}=(0,-1)$ and $v^{}_{2}(s)=-1+2s$, with $0 \leq s \leq 1$. The second term in the RHS of~\eqref{WFh-eq-11} yields
\begin{align*}
& -\int^{}_{\Gamma_{1}^{}} \Big(\frac{1}{2} \big( u^{}_{\,\Gamma_{1}^{}}|^{}_{K^{}_{1}} + u^{}_{\,\Gamma_{1}^{}}|^{}_{K^{}_{3}} \big)  \mathbf{a} \cdot \mathbf{n}_{1}^{} 
+ \frac{1}{2} \mid \! \mathbf{a} \cdot \mathbf{n}_{1}^{} \! \mid \big( u^{}_{\,\Gamma_{1}^{}}|^{}_{K^{}_{1}} - u^{}_{\,\Gamma_{1}^{}}|^{}_{K^{}_{3}} \big) \Big) v^{}_{2}(s) ds \\
& = -\int^{}_{\Gamma_{1}^{}} \Big(-\frac{a_{2}^{}}{2} \big( u^{}_{\,\Gamma_{1}^{}}|^{}_{K^{}_{1}} + u^{}_{\,\Gamma_{1}^{}}|^{}_{K^{}_{3}} \big) + \frac{a_{2}^{}}{2} \big( u^{}_{\,\Gamma_{1}^{}}|^{}_{K^{}_{1}} - u^{}_{\,\Gamma_{1}^{}}|^{}_{K^{}_{3}} \big) \Big) v^{}_{2}(s) ds \\
&=a_{2}^{} \int^{}_{\Gamma_{1}^{}} u^{}_{\,\Gamma_{1}^{}}|^{}_{K^{}_{3}} v^{}_{2}(s) ds =  a_{2}^{} h \int^{1}_{0} u^{}_{\,\Gamma_{1}^{}}|^{}_{K^{}_{3}} v^{}_{2}(s) ds = \frac{a_{2}^{} h}{3} (- u^{}_{15} + u^{}_{17}), 
\label{comp-eq-NC-node-2-gam-1-transp}
\end{align*}
since $u^{}_{\Gamma_{1}^{}}|^{}_{K^{}_{3}} = u^{}_{1} + u^{}_{15} - u^{}_{17} - 2s ( u^{}_{15} - u^{}_{17} )$, with $0 \leq s \leq 1$.  

Following the same procedure along $\Gamma_{1}^{}$ but for the last term in the RHS of~\eqref{WFh-eq-11}, with $p/h^{}_{ed}=p/h^{}_{\,\Gamma_{1}^{}}=p/h$ and $\mathbf{n}_{1}^{}=(0,-1)$, yielding $\sign (\mathbf{b} \cdot \mathbf{n}^{}_{1}) < 0$, we obtain
\begin{align*}
& \int^{}_{\Gamma_{1}^{}} \!\!\Big(\frac{1}{2} \big( \boldsymbol{\theta}^{}_{\Gamma_{1}^{}}|^{}_{K^{}_{1}}\! + \boldsymbol{\theta}^{}_{\Gamma_{1}^{}}|^{}_{K^{}_{3}} \big) + \frac{1}{2} \big( \boldsymbol{\theta}^{}_{\Gamma_{1}^{}}|^{}_{K^{}_{1}}\! - \boldsymbol{\theta}^{}_{\Gamma_{1}^{}}|^{}_{K^{}_{3}} \big) - \frac{p}{h} \big( u^{}_{\,\Gamma_{1}^{}}|^{}_{K^{}_{1}} \!\!- u^{}_{\,\Gamma_{1}^{}}|^{}_{K^{}_{3}} \big) \mathbf{n}^{}_{1} \Big) \! \cdot \mathbf{n}^{}_{1} v^{}_{2}(s) ds \\
& = \int^{}_{\Gamma_{1}^{}} \Big(\boldsymbol{\theta}^{}_{\Gamma_{1}^{}}|^{}_{K^{}_{1}} + \frac{p}{h} \big( u^{}_{\,\Gamma_{1}^{}}|^{}_{K^{}_{3}} - u^{}_{\,\Gamma_{1}^{}}|^{}_{K^{}_{1}} \big) \mathbf{n}^{}_{1} \Big) \cdot \mathbf{n}^{}_{1} v^{}_{2}(s) ds \\ 
& =  -h \int^{1}_{0} \Big( \theta^{y}_{1} - \theta^{y}_{2} + \theta^{y}_{3} - 2s \big( \theta^{y}_{3} - \theta^{y}_{2} \big) \Big) (2s-1) ds \\
& + p \!\int^{1}_{0}\! \Big( \big( u^{}_{1} + u^{}_{15} - u^{}_{17} - 2s (u ^{}_{15} - u^{}_{17} ) \big) - \big(u ^{}_{1} - u^{}_{2} + u^{}_{3} - 2s ( u^{}_{3} - u^{}_{2} ) \big) \Big) (2s-1) ds \\
& = \frac{h}{3} \big( - \theta^{y}_{2} + \theta^{y}_{3} \big) + \frac{p}{3} \big( - u^{}_{2} + u^{}_{3} - u^{}_{15} + u^{}_{17} \big).
\label{comp-eq-NC-node-2-gam-1}
\end{align*}
The computations of the line integrals in the RHS of~\eqref{WFh-eq-11} (second and last terms) are conducted similarly on $\Gamma_{3}^{}$, $\Gamma_{4}^{}$ and $\Gamma_{5}^{}$. 

The equation for $\boldsymbol{\theta}_{h}^{}$ in~\eqref{WFh-eq-2} is also computed at node $2$ of~\cref{Node-distrib-P1NC} leading to
\begin{align}
\label{WFh-eq-22}
\int^{}_{K^{}_{1} \cup K^{}_{2}} \boldsymbol{\theta}_{h}^{} \cdot \boldsymbol{\psi}^{}_{2} \, d\mathbf{x} = - \nu \int^{}_{K^{}_{1} \cup K^{}_{2}} u_{h}^{} \, \nabla \cdot \boldsymbol{\psi}^{}_{2}\, d\mathbf{x} + \nu \int^{}_{\partial K^{}_{1} \cup K^{}_{2}} u_{h}^{*} \, \boldsymbol{\psi}^{}_{2} \cdot \mathbf{n}^{}_{ed} \, ds.
\end{align}
By taking advantage of the previous calculations we obtain
\begin{align*}
\int^{}_{K^{}_{1} \cup K^{}_{2}} \boldsymbol{\theta}_{h}^{} \cdot \boldsymbol{\psi}^{}_{2} \,d{\mathbf{x}}
& = \frac{h_{}^{2}}{3} {\theta}_{2}^{x} \mbox{ for } \boldsymbol{\psi}^{}_{2}=({\psi}^{}_{2},0) \mbox{ and } \frac{h_{}^{2}}{3} {\theta}_{2}^{y} \mbox{ for } \boldsymbol{\psi}^{}_{2}=(0,{\psi}^{}_{2}), \\
- \int^{}_{K^{}_{1} \cup K^{}_{2}} u_{h}^{} \, \nabla \cdot \boldsymbol{\psi}^{}_{2}\, d\mathbf{x}
& = \frac{h}{3} (- u^{}_{1}- u^{}_{3} + u^{}_{4} + u^{}_{5}) \mbox{ for } \boldsymbol{\psi}^{}_{2}=({\psi}^{}_{2},0)\mbox{ and } \boldsymbol{\psi}^{}_{2}=(0,{\psi}^{}_{2}). 
\end{align*}
In the case $\boldsymbol{\psi}^{}_{2}=({\psi}^{}_{2},0)$, and along $\Gamma_{3}^{}$ for example, with $\mathbf{n}_{3}^{}=(-1,0)$, $\sign (\mathbf{b} \cdot \mathbf{n}^{}_{3}) < 0$ and $\psi^{}_{2}(s)=1-2s$ for $0 \leq s \leq 1$, the last term in the RHS of~~\eqref{WFh-eq-22} yields
\begin{align*}
& \int^{}_{\Gamma_{3}^{}} \!\!\Big(\frac{1}{2} \big( u^{}_{\,\Gamma_{3}^{}}|^{}_{K^{}_{1}}\! + u^{}_{\,\Gamma_{3}^{}}|^{}_{K^{}_{4}} \big) 
+ \frac{1}{2} \big( u^{}_{\,\Gamma_{3}^{}}|^{}_{K^{}_{4}} - u^{}_{\,\Gamma_{3}^{}}|^{}_{K^{}_{1}} \big) 
                        \Big) \boldsymbol{\psi}^{}_{2} \cdot \mathbf{n}^{}_{3} ds \\
& =  h \int^{1}_{0} \Big( u^{}_{3} + u^{}_{10} - u^{}_{11} - 2s \big( u^{}_{10} - u^{}_{11} \big) \Big) (2s-1) ds  =\frac{h}{3} \big(-u^{}_{10} + u^{}_{11} \big).
\label{comp-eq-NC-node-2-gam-3}
\end{align*}
The computations of the line integrals in the RHS of~\eqref{WFh-eq-22} are conducted similarly on $\Gamma_{1}^{}$, $\Gamma_{4}^{}$ and $\Gamma_{5}^{}$. These are null on $\Gamma_{1}^{}$ and $\Gamma_{5}^{}$ in the case $\boldsymbol{\psi}^{}_{2}=({\psi}^{}_{2},0)$ because the first component of $\mathbf{n}^{}_{1}$ and $\mathbf{n}^{}_{5}$ is equal to zero. Since the numerical trace $u_{h}^{*}$ in~\eqref{jump_u} depends only on $u_{h}^{}$, the values of $\theta_{h}^{x}$ and $\theta_{h}^{y}$ also depend only on $u_{h}^{}$. That is why the method is referred to as the \textit{local} NC method, named LNC in this paper.

In order to perform the Fourier analysis, the discrete equations for $u_{h}^{}$ in~\eqref{WFh-eq-11} and $\boldsymbol{\theta}_{h}^{}$ in~\eqref{WFh-eq-22} are computed similarly at nodes $1$ and $3$ of~\cref{Node-distrib-P1NC}. The vectors $\boldsymbol{\Lambda}^{}_{1}, \boldsymbol{\Lambda}^{}_{2}, \boldsymbol{\Lambda}^{}_{3}, \boldsymbol{\Lambda}^{}_{4}$ in~\cref{compact-discrete-eqn-u,compact-discrete-eqn-qy} are shown in~\cref{app-def-submatNC}.
\subsection{The \boldmath$P_{1}^{DG}$ discrete equations}
\label{DG-discrete}
\begin{figure}[h!]
\begin{center}
\begin{tabular}{ccc}
\includegraphics*[height=28mm]{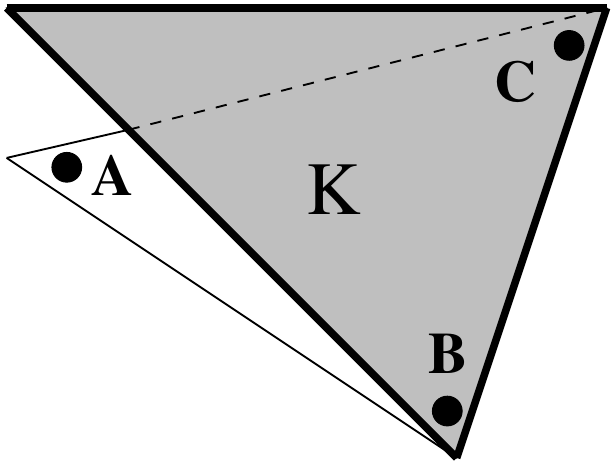}
&
\includegraphics*[height=28mm]{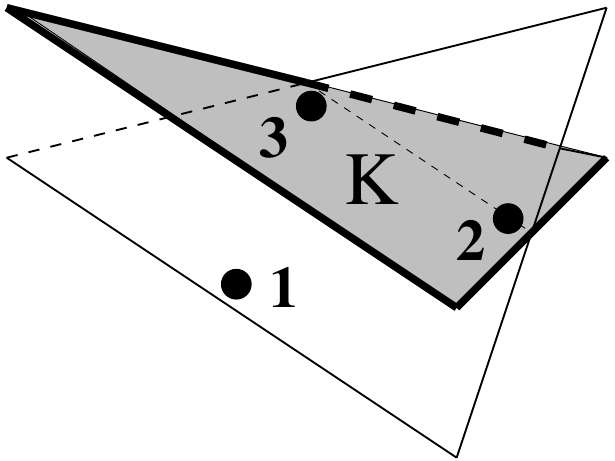}
&
\includegraphics*[height=42mm]{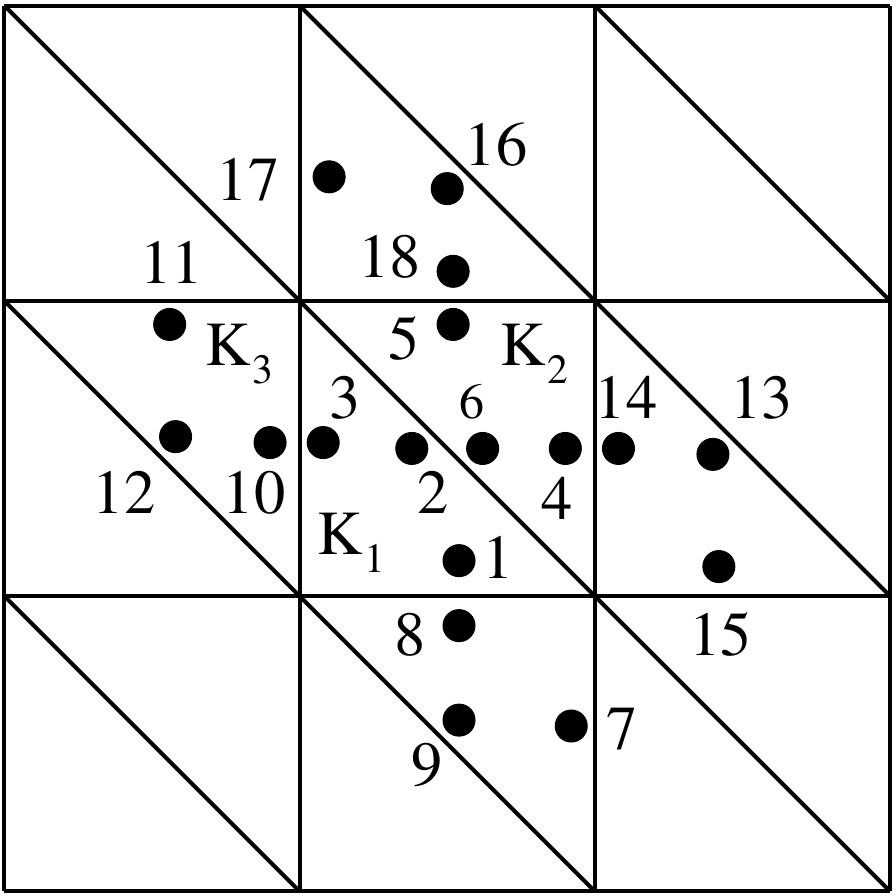}
\end{tabular}
\caption{The $P_{1}^{}$ conforming linear basis function $\varphi$ at node $A$ in grey (left). The linear basis function $v$ or $\psi$ at node $1$ in grey (middle). Node locations, used to compute the discrete equations, are represented by the symbol {\Large $\bullet$} (right).}
\label{Node-distrib-P1DG}
\end{center}
\end{figure}

The space $V_{h,K_{el}^{}}^{} \subset V_{K_{el}^{}}^{}$ is now defined~as
\begin{equation*}
V_{h,K_{el}^{}}^{} := V_{h,K_{el}^{}}^{DG} = \{u_{h}^{} : u_{h}^{} \vert^{}_{K_{el}^{}} \in \mathcal{P}_{}^{1} (K_{el}^{}) \mbox{ for all } K_{el}^{} \in \tau_{h}^{}\}.
\end{equation*}
The discrete solutions $u_{h}^{}$ and $\boldsymbol{\theta}_{h}^{}$ are again approximated linearly inside each triangle $K_{el}^{}$ of $\tau_{h}^{}$. However, the nodal values $u_{j}^{}$ and $\boldsymbol{\theta}_{j}^{}$, $j=1,2,3,...$, now exclusively belong to a given element $K_{el}^{}$ and hence such an approximation is completely discontinuous between the triangles of $\tau_{h}^{}$. There are several choices for the location of the nodal values \textit{inside} $K_{el}^{}$, and hence for the basis functions. One possibility is to locate the unknown at the vertices of $K_{el}^{}$ using the basis function $\varphi$, mentioned in~\cref{NC-discrete}, and displayed in~\cref{Node-distrib-P1DG} (left) at vertex $A$. A second choice is to locate the unknown at the mid-edge points of $K_{el}^{}$ using the basis functions $v$ and $\psi$, shown in~\cref{Node-distrib-P1DG} (middle) at mid-edge point $1$. But contrary to the NC discretisation performed in~\cref{NC-discrete}, there are now three mid-edge points per element $K_{el}^{}$, and hence $u_{h}^{}$ and $\boldsymbol{\theta}_{h}^{}$ are no longer continuous at the mid-edges of $K_{el}^{}$. Both choices lead to the same dispersion relation obtained in~\cref{Four-disp}. The second one is adopted in this study as it yields to diagonal mass matrices and hence it simplifies the writing of the discrete equations.

The equation for $u_{h}^{}$ in~\eqref{WFh-eq-1} is computed at node $2$ of~\cref{Node-distrib-P1DG} (right panel), and the support of the basis function $\psi^{}_{2}$ is now only made up of triangle $K^{}_{1}$, yielding
\begin{align}
\label{WFh-eq-31}
\int^{}_{K^{}_{1}} \frac{\partial u_{h}^{}}{\partial t} \, v^{}_{2} \, d\mathbf{x} & = \int^{}_{K^{}_{1}} u_{h}^{}\, \mathbf{a} \cdot \nabla v^{}_{2} \, d\mathbf{x} - \int^{}_{\partial K^{}_{1}} (u_{h}\, \mathbf{a})_{}^{\star} \cdot \mathbf{n}^{}_{ed}\,  v^{}_{2} \, ds \\ 
& - \int^{}_{K^{}_{1}} \boldsymbol{\theta}_{h}^{} \cdot 
\nabla v^{}_{2} \, d\mathbf{x} + \int^{}_{\partial K^{}_{1}}  \boldsymbol{\theta}_{h}^{*} \cdot \mathbf{n}^{}_{ed}\, v^{}_{2} \, ds. \nonumber 
\end{align}
The integrals in~\eqref{WFh-eq-31} are computed by expanding $u_{h}^{}$ and $\boldsymbol{\theta}_{h}^{}$ over $K_{1}^{}$ in terms of the basis functions $v^{}_{1}, v^{}_{2}, v^{}_{3}$ and $\psi^{}_{1}, \psi^{}_{2}, \psi^{}_{3}$. Since the bases $\{v^{}_{j}\}$ and $\{\psi^{}_{j}\}$, $j=1,2,3,...$, are orthogonal sets and elementary computations yield $\frac{\partial v^{}_{2}}{\partial x}|^{}_{K_1} = \frac{\partial v^{}_{2}}{\partial y}|^{}_{K_1} = \frac{2}{h}$, we obtain
\begin{align*}
\int^{}_{K^{}_{1}} \frac{\partial u_{h}^{}}{\partial t}\,v^{}_{2}\,d{\mathbf{x}}
& = \frac{h_{}^{2}}{6} \frac{\partial u_{2}^{}}{\partial t}, \\
\int^{}_{K^{}_{1}} u_{h}^{}\, \mathbf{a} \cdot \nabla v^{}_{2} \, d\mathbf{x} 
& = \frac{h}{3} (a_{1}^{}+a_{2}^{}) (u_{1}^{}+u_{2}^{}+u_{3}^{}) \\
- \int^{}_{K^{}_{1}} \boldsymbol{\theta}_{h}^{} \cdot \nabla v^{}_{2} \, d\mathbf{x} 
& = - \frac{h}{3} \big( \theta^{x}_{1} + \theta^{x}_{2} + \theta^{x}_{3} + \theta^{y}_{1} + \theta^{y}_{2} +\theta ^{y}_{3} \big).
\end{align*}
We have $\partial K^{}_{1} = \Gamma_{1}^{} \cup \Gamma_{2}^{}\cup \Gamma_{3}^{}$ and along $\Gamma_{2}^{}$ for example, $\mathbf{n}_{2}^{}=\frac{1}{\sqrt{2}} (1,1)$ and $\psi^{}_{2}(s)=1$ for $0 \leq s \leq 1$. The second term in the RHS of~\eqref{WFh-eq-31} yields
\begin{align*}
& -\int^{}_{\Gamma_{2}^{}} \Big(\frac{1}{2} \big( u^{}_{\,\Gamma_{2}^{}}|^{}_{K^{}_{1}} + u^{}_{\,\Gamma_{2}^{}}|^{}_{K^{}_{2}} \big)  \mathbf{a} \cdot \mathbf{n}_{2}^{} + \frac{1}{2} \mid \! \mathbf{a} \cdot \mathbf{n}_{2}^{} \! \mid \big( u^{}_{\,\Gamma_{2}^{}}|^{}_{K^{}_{1}} - u^{}_{\,\Gamma_{2}^{}}|^{}_{K^{}_{2}} \big) \Big) v^{}_{2}(s) ds \\
& = -\int^{}_{\Gamma_{2}^{}} \Big(\frac{(a_{1}^{}+a_{2}^{})}{2\sqrt{2}} \big( u^{}_{\,\Gamma_{2}^{}}|^{}_{K^{}_{1}} + u^{}_{\,\Gamma_{2}^{}}|^{}_{K^{}_{3}} \big) + \frac{(a_{1}^{}+a_{2}^{})}{2\sqrt{2}} \big( u^{}_{\,\Gamma_{2}^{}}|^{}_{K^{}_{1}} - u^{}_{\,\Gamma_{2}^{}}|^{}_{K^{}_{3}} \big) \Big) v^{}_{2}(s) ds \\
&=- \frac{(a_{1}^{}+a_{2}^{})}{\sqrt{2}} \int^{}_{\Gamma_{2}^{}} u^{}_{\,\Gamma_{2}^{}}|^{}_{K^{}_{1}} v^{}_{2}(s) ds =  -(a_{1}^{}+a_{2}^{}) h \int^{1}_{0} u^{}_{\,\Gamma_{2}^{}}|^{}_{K^{}_{1}} v^{}_{2}(s) ds = -(u+v) h u_{2}^{}, 
\label{comp-eq-DG-node-2-gam-2-transp}
\end{align*}
since $u^{}_{\,\Gamma_{2}^{}}|^{}_{K^{}_{1}} = u^{}_{2} + \big(u^{}_{1}- u^{}_{3} \big) (1-2s)$, with $0 \leq s \leq 1$.  

Following the same procedure along $\Gamma_{2}^{}$ but for the last term in the RHS of~\eqref{WFh-eq-31}, with $v^{}_{2}(s)=1$, $p/h^{}_{ed}=p/h^{}_{\,\Gamma_{2}^{}}=p/(h\sqrt{2})$ and $\sign (\mathbf{b} \cdot \mathbf{n}^{}_{2}) > 0$, we obtain
\begin{align*}
& \int^{}_{\Gamma_{2}^{}} \!\!\Big(\frac{1}{2} \big( \boldsymbol{\theta}^{}_{\Gamma_{2}^{}}|^{}_{K^{}_{1}}\! + \boldsymbol{\theta}^{}_{\Gamma_{6}^{}}|^{}_{K^{}_{2}} \big) - \frac{1}{2} \big( \boldsymbol{\theta}^{}_{\Gamma_{2}^{}}|^{}_{K^{}_{1}}\! - \boldsymbol{\theta}^{}_{\Gamma_{6}^{}}|^{}_{K^{}_{2}} \big) - \frac{p}{h\sqrt{2}} \big( u^{}_{\,\Gamma_{2}^{}}|^{}_{K^{}_{1}} \!- u^{}_{\,\Gamma_{6}^{}}|^{}_{K^{}_{2}} \big) \mathbf{n}^{}_{2} \Big) \! \cdot \mathbf{n}^{}_{2} ds \\
& = \int^{}_{\Gamma_{2}^{}} \Big(\boldsymbol{\theta}^{}_{\Gamma_{6}^{}}|^{}_{K^{}_{2}} + \frac{p}{h\sqrt{2}} \big( u^{}_{\,\Gamma_{6}^{}}|^{}_{K^{}_{2}} - u^{}_{\,\Gamma_{2}^{}}|^{}_{K^{}_{1}} \big) \mathbf{n}^{}_{2} \Big) \cdot \mathbf{n}^{}_{2} ds \\ 
& =  h \int^{1}_{0} \Big( \theta^{x}_{6} + \big(\theta^{x}_{4} - \theta^{x}_{5}\big) (1-2s) + \theta^{y}_{6} + \big(\theta^{y}_{4} - \theta^{y}_{5}\big) (1-2s) \Big) ds \\
& \quad + p \int^{1}_{0} \Big( u^{}_{6} + \big(u^{}_{4}- u^{}_{5} \big) (1-2s) - u^{}_{2} - \big(u^{}_{1}- u^{}_{3} \big) (1-2s) \Big) ds \\*[0.1cm]
& = h \big( \theta^{x}_{6} + \theta^{y}_{6} \big) + p \big(u^{}_{6} - u^{}_{2} \big),
\label{comp-eq-NC-local-node-2-gam-1}
\end{align*}
since for example, $u^{}_{\,\Gamma_{6}^{}}|^{}_{K^{}_{2}} = u^{}_{6} + \big(u^{}_{4}- u^{}_{5} \big) (1-2s)$, with $0 \leq s \leq 1$. The computations of the two line integrals in the RHS of~\eqref{WFh-eq-31} are conducted similarly on $\Gamma_{1}^{}$ and $\Gamma_{3}^{}$.

The equation for $\boldsymbol{\theta}_{h}^{}$ in~\eqref{WFh-eq-2} is also computed at node $2$ of~\cref{Node-distrib-P1DG} leading to
\begin{align}
\label{WFh-eq-222}
\int^{}_{K^{}_{1}} \boldsymbol{\theta}_{h}^{} \cdot \boldsymbol{\psi}^{}_{2} \, d\mathbf{x} = - \nu \int^{}_{K^{}_{1}} u_{h}^{} \, \nabla \cdot \boldsymbol{\psi}^{}_{2}\, d\mathbf{x} + \nu \int^{}_{\partial K^{}_{1}} u_{h}^{*} \, \boldsymbol{\psi}^{}_{2} \cdot \mathbf{n}^{}_{ed} \, ds.
\end{align}
By taking advantage of the previous calculations for $u_{h}^{}$ we obtain
\begin{align*}
\int^{}_{K^{}_{1}} \boldsymbol{\theta}_{h}^{} \cdot \boldsymbol{\psi}^{}_{2} \,d{\mathbf{x}}
& = \frac{h_{}^{2}}{6} {\theta}_{2}^{x} \mbox{ for } \boldsymbol{\psi}^{}_{2}=({\psi}^{}_{2},0) \mbox{ and } \frac{h_{}^{2}}{6} {\theta}_{2}^{y } \mbox{ for } \boldsymbol{\psi}^{}_{2}=(0,{\psi}^{}_{2}), \\
- \int^{}_{K^{}_{1}} u_{h}^{} \frac{\partial {\psi}^{}_{2}}{\partial x}\, d\mathbf{x} = &- \int^{}_{K^{}_{1}} u_{h}^{} \frac{\partial {\psi}^{}_{2}}{\partial y}\, d\mathbf{x}
 = -\frac{h}{3} (u^{}_{1}+ u^{}_{2} + u^{}_{3}). 
\end{align*}
In the case $\boldsymbol{\psi}^{}_{2}=({\psi}^{}_{2},0)$, and along $\Gamma_{3}^{}$ for example, $\mathbf{n}_{3}^{}=(-1,0)$, $\sign (\mathbf{b} \cdot \mathbf{n}^{}_{3}) < 0$ and $\psi^{}_{2}(s)=1-2s$ for $0 \leq s \leq 1$. The last term in the RHS of~~\eqref{WFh-eq-222} yields
\begin{align*}
& \int^{}_{\,\Gamma_{3}^{}}\Big(\frac{1}{2} \big( u^{}_{\,\Gamma_{3}^{}}|^{}_{K^{}_{1}} + u^{}_{\,\Gamma_{3}^{}}|^{}_{K^{}_{3}} \big) 
+ \frac{1}{2} \big( u^{}_{\,\Gamma_{3}^{}}|^{}_{K^{}_{3}} - u^{}_{\Gamma_{3}^{}}|^{}_{K^{}_{1}} \big) 
                        \Big) \boldsymbol{\psi}^{}_{2} \cdot \mathbf{n}^{}_{3} ds \\
& =  h \int^{1}_{0} \Big( u^{}_{10} + u^{}_{11} - u^{}_{12}- 2s \big( u^{}_{11} - u^{}_{12} \big) \Big) (2s-1) ds  =\frac{h}{3} \big(-u^{}_{11} + u^{}_{12} \big),
\label{comp-eq-DG-node-2-gam-3}
\end{align*}
since $u^{}_{\,\Gamma_{3}^{}}|^{}_{K^{}_{3}}=u^{}_{10} + u^{}_{11} - u^{}_{12}- 2s \big( u^{}_{11} - u^{}_{12})$. The computations of the line integrals in the RHS of~\eqref{WFh-eq-222} are conducted similarly on $\Gamma_{1}^{}$ and $\Gamma_{3}^{}$.

In order to perform the Fourier analysis, the discrete equations for $u_{h}^{}$ in~\eqref{WFh-eq-31} and $\boldsymbol{\theta}_{h}^{}$ in~\eqref{WFh-eq-222} are also computed at nodes $1, 3, 4, 5$ and $6$ of~\cref{Node-distrib-P1DG}. The vectors  $\boldsymbol{\Lambda}^{}_{1}, \boldsymbol{\Lambda}^{}_{2}, \boldsymbol{\Lambda}^{}_{3}, \boldsymbol{\Lambda}^{}_{4}$ in~\cref{compact-discrete-eqn-u,compact-discrete-eqn-qy} are shown in~\cref{app-def-submatDG}.
\section{The Fourier / dispersion analysis}
\label{Four-disp}
\subsection{The generalized eigenvalue system}
The Fourier / dispersion analysis has proven practical and beneficial to detect dispersion and dissipation problems as well as the presence of spurious modes in a few discretization schemes \cite{Cotter-Ham-2011, LeRoux-SIAM-2005, LeRoux-SIAM-2007, LeRoux-JCP-2012, LeRoux-JCP-2024}. It is adopted in this study to analyse the stability and accuracy of the $P^{NC}_{1}$ and $P^{DG}_{1}$ schemes examined earlier in~\cref{numerical:heat}. Note that the mesh length parameter $h$ remains constant in the present section.

\begin{figure}[h!]
\begin{center}
{\setlength{\tabcolsep}{10pt}
\begin{tabular}{cc|ccc}
\multicolumn{2}{c}{$P^{NC}_{1}$} & \multicolumn{3}{c}{$P^{DG}_{1}$} \\*[0.2cm] 
\vrule height 10mm width 0 pt depth 0mm
\begin{minipage}{17mm}
\begin{center}
\includegraphics*[width=17mm]{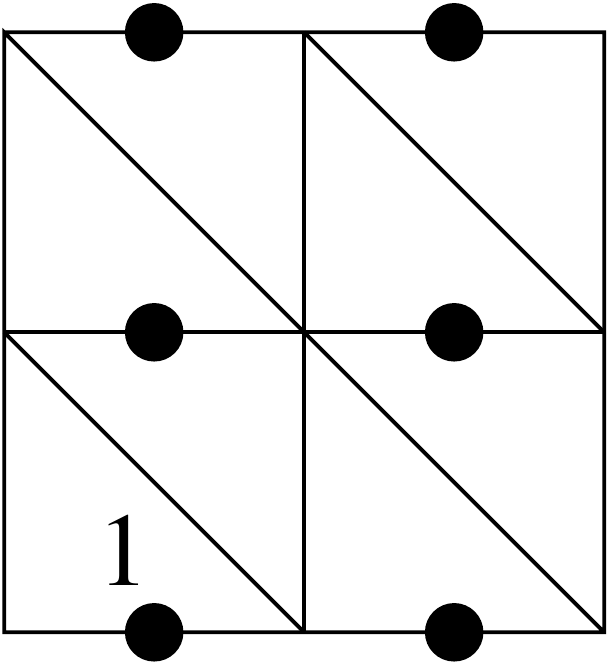}
\end{center}
\end{minipage}
&
\begin{minipage}{17mm}
\begin{center}
\includegraphics*[width=17mm]{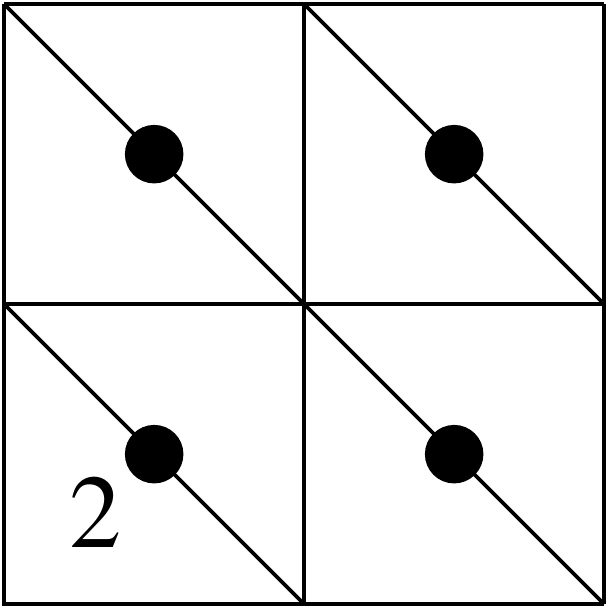}
\end{center}
\end{minipage}
&
\begin{minipage}{17mm}
\begin{center}
\includegraphics*[width=17mm]{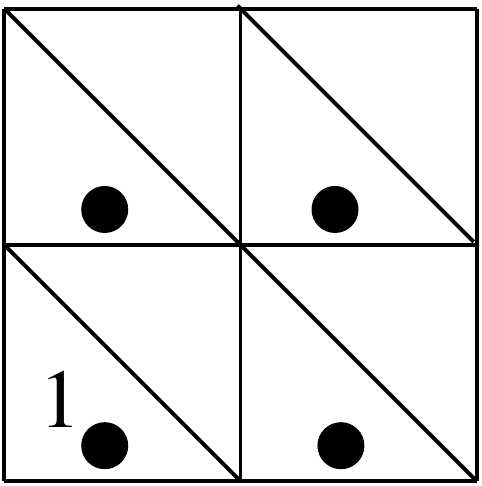}
\end{center}
\end{minipage}
&
\begin{minipage}{17mm}
\begin{center}
\includegraphics*[width=17mm]{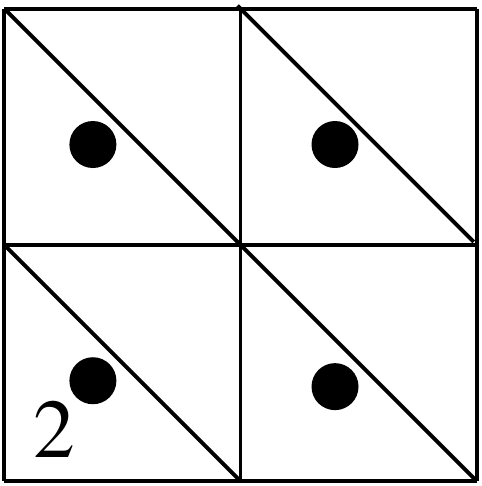}
\end{center}
\end{minipage}
&
\begin{minipage}{17mm}
\begin{center}
\includegraphics*[width=17mm]{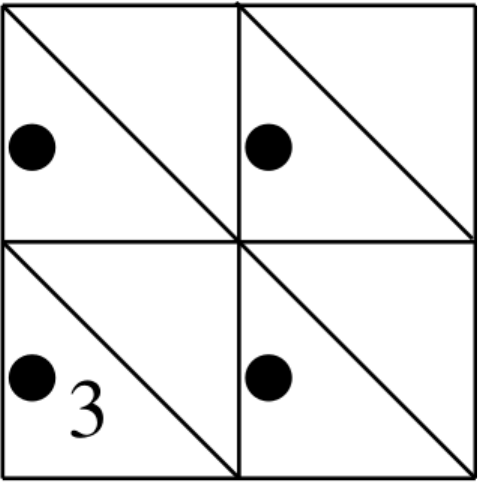}
\end{center}
\end{minipage}
\\*[0.8cm] 
\begin{minipage}{18mm}
\begin{center}
\includegraphics*[width=18mm]{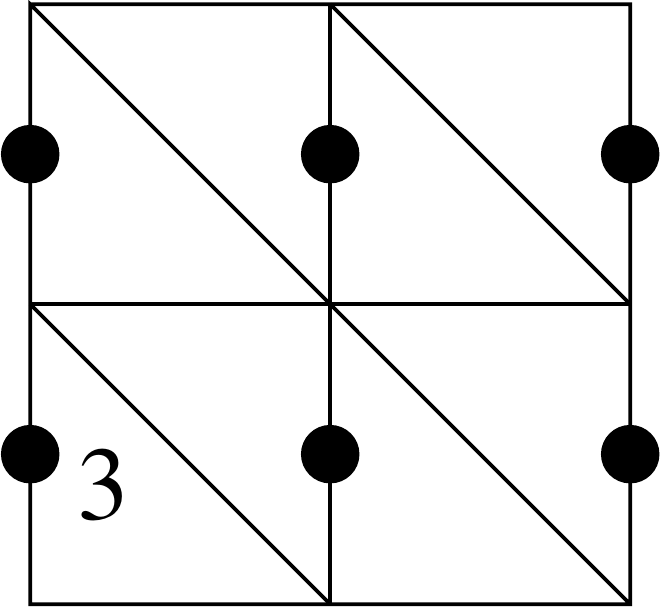}
\end{center}
\end{minipage}
&
&
\begin{minipage}{17mm}
\begin{center}
\includegraphics*[width=17mm]{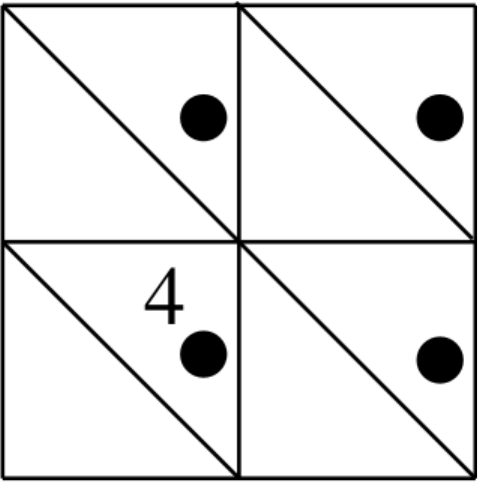}
\end{center}
\end{minipage}
&
\begin{minipage}{17mm}
\begin{center}
\includegraphics*[width=17mm]{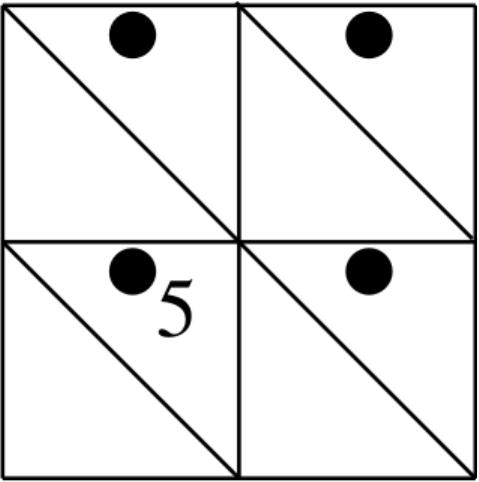}
\end{center}
\end{minipage}
&
\begin{minipage}{17mm}
\begin{center}
\includegraphics*[width=17mm]{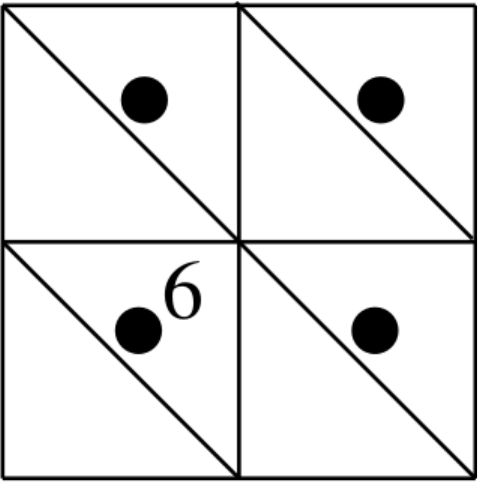}
\end{center}
\end{minipage}
\end{tabular}
}
\caption{The sets of nodal unknowns for the $P^{NC}_{1}$ and $P^{DG}_{1}$ schemes. The nodal unknowns belonging to the same set are distributed on a regular grid of size $h$.}
\label{Node-distrib-Fourier}
\end{center}
\end{figure}

When continuous linear finite elements are employed, the Fourier analysis is straightforward since all the discrete unknowns are located at only one type of mesh nodes, namely vertices. It is hence sufficient to consider one discrete equation for $u_{h}^{}$ and $\boldsymbol{\theta}_{h}^{}$ at such a typical node to compute the dispersion relation, due to symmetry reasons. In this study, such a simple case does not occur since the nodal unknowns appearing in~\cref{compact-discrete-eqn-u,compact-discrete-eqn-qx,compact-discrete-eqn-qy} are located at midpoints of edges for the $P^{NC}_{1}$ method and local midpoints of edges per triangle for the $P^{DG}_{1}$ scheme. Consequently, we need to compute the discrete equations for $u_{h}^{}$ and $\boldsymbol{\theta}_{h}^{}$ at the three types of faces: namely, horizontal, vertical and diagonal faces for the $P^{NC}_{1}$ scheme and at the six types of midpoints of edges for the $P^{DG}_{1}$ scheme, namely three local midpoints of edges per triangle for the upward and downward pointing triangles. This explains why, equations~\cref{compact-discrete-eqn-u,compact-discrete-eqn-qx,compact-discrete-eqn-qy} have been computed at nodes $1,2,3,$ of~\cref{Node-distrib-P1NC} in~\cref{NC-discrete} for the $P^{NC}_{1}$ method and at nodes $1,2,3,4,5,6,$ in~\cref{Node-distrib-P1DG} for the $P^{DG}_{1}$ scheme. Note that the nodal unknowns belonging to the same set are distributed on a regular grid of size $h$, as shown in~\cref{Node-distrib-Fourier}. 
\setlength{\arrayrulewidth}{0.7pt}

In order to compute the dispersion relation in~\cref{Disp-Rel} at the discrete level, we follow the same procedure as in~\cref{Problem-model} for the continuous case. The discrete solutions corresponding to 
\begin{equation}
{u}^{}_{j} = \widehat{u}^{}_{j^{}_0}\,e^{i(kx_{j}^{}+ly_{j}^{} + \omega t)}, \quad \theta_{j}^{x}=\widehat{\theta}^{\,x}_{j^{}_0}\,e^{i(kx_{j}^{}+ly_{j}^{} + \omega t)}, \quad \theta_{j}^{y}=\widehat{\theta}^{\,y}_{j^{}_0}\,e^{i(kx_{j}^{}+ly_{j}^{} + \omega t)},
\label{Four-expansion}
\end{equation}
are sought, where $u_{j}^{}, \theta_{j}^{x}$ and $\theta_{j}^{y}$ are the nodal unknowns that appear in the selected discrete equations~\cref{compact-discrete-eqn-u,compact-discrete-eqn-qx,compact-discrete-eqn-qy}, and $\widehat{u}^{}_{j^{}_0}, \widehat{\theta}^{\,x}_{j^{}_0}, \widehat{\theta}^{\,y}_{j^{}_0}$, for $j^{}_0 \in \mathcal{J}_{n}^{}:=\left \{1,2,...,n \right \}$, are the Fourier amplitudes at the typical nodes $1,2,3,$ of~\cref{Node-distrib-P1NC} for the $P^{NC}_{1}$ method ($n=3$) and at nodes $1,2,3,4,5,6,$ of~\cref{Node-distrib-P1DG} for the $P^{DG}_{1}$ scheme ($n=6$). Further, $k$ and $l$ are the wave numbers in the $x$- and $y$-directions, respectively, $\omega$ is the harmonic frequency, and the $(x_{j}^{},y_{j}^{})$ coordinates are expressed in terms of a distance to a reference node. 

For example, considering~\eqref{compact-discrete-eqn-qx} at node $1$ by using $\boldsymbol{\Lambda}^{}_{1}$ in~\cref{app-def-submatNC} yields
\begin{equation}
\theta_{1}^{x} = \frac{\nu}{h} ( 2 u^{}_{1} -   u^{}_{3} + u^{}_{10} - u^{}_{11} - u^{}_{14} - u^{}_{15} + u^{}_{16}).
\label{q1discrete}
\end{equation}
We let $\alpha = e_{}^{\frac{ikh}{2}}$ and $\beta = e_{}^{\frac{ilh}{2}}$. The Fourier expansions~\eqref{Four-expansion} are then inserted in~\eqref{q1discrete} leading to
\begin{align*}
\widehat{\theta}^{\,x}_{1} = \frac{\nu}{h} & \left( 2 \widehat{u}^{}_{1} - \widehat{u}^{}_{3} e_{}^{i(-k+l)\frac{h}{2}}   
+ \widehat{u}^{}_{1} e_{}^{i(-k+l)h} - \widehat{u}^{}_{2} e_{}^{i(-2k+l)\frac{h}{2}} \right. \\
& \left. - \widehat{u}^{}_{3} e_{}^{i(-k-l)\frac{h}{2}} - \widehat{u}^{}_{2} e_{}^{-il\frac{h}{2}} + \widehat{u}^{}_{1} e_{}^{-ilh} \right), \\
= \frac{\nu}{h} & \Big( 2 \widehat{u}^{}_{1} - \widehat{u}^{}_{3}\frac{\beta}{\alpha} + \widehat{u}^{}_{1} \frac{\beta^{2}_{}}{\alpha^{2}_{}} - \widehat{u}^{}_{2} \frac{\beta}{\alpha^{2}_{}}- \widehat{u}^{}_{3} \frac{1}{\alpha \beta} - \widehat{u}^{}_{2} \frac{1}{\beta} + \widehat{u}^{}_{1} \frac{1}{\beta^{2}_{}} \Big),
\end{align*}
since the amplitudes at nodes $10$ and $16$ (resp. at nodes $11$ and $15$) behave in the same manner as the amplitude at node $1$ (resp. node $2$), and the amplitude at node $14$ behaves as for the amplitude at node $3$, for symmetry reason. 

By following the same procedure, substitution of $(u_{j}^{}, \theta_{j}^{x},\theta_{j}^{y})$ from~\eqref{Four-expansion} into the $3n$ selected discrete equations~\cref{compact-discrete-eqn-u,compact-discrete-eqn-qx,compact-discrete-eqn-qy}, for the $P^{NC}_{1}$ (with $n=3$) and $P^{DG}_{1}$ (with $n=6$) schemes, leads to the $3n \times 3n$ generalized eigenvalue problem  
\begin{equation}
{A}_{3n,3n}^{}\,\widehat{\mathbf{X}} = i \omega\,{B}_{3n,3n}^{}\,\widehat{\mathbf{X}}, 
\label{disprel}
\end{equation}
with $\widehat{\mathbf{X}} = \left(
\{\widehat{u}^{}_{j^{}_0}\}^{}_{j^{}_0 \in \mathcal{J}_{n}^{}}, \{\widehat{\theta}^{\,x}_{j^{}_0}\}^{}_{j^{}_0 \in \mathcal{J}_{n}^{}}, \{\widehat{\theta}^{\,y}_{j^{}_0}\}^{}_{j^{}_0 \in \mathcal{J}_{n}^{}}
\right )$, \quad $
{B}_{3n,3n}^{} =
\begin{pmatrix} 
I_{n,n}^{} & O_{n,n}^{} & O_{n,n}^{} \\*[0.1cm]
O_{n,n}^{} & O_{n,n}^{} & O_{n,n}^{} \\*[0.1cm]
O_{n,n}^{} & O_{n,n}^{} & O_{n,n}^{}
\end{pmatrix}, 
$
\begin{equation}
\label{matrix-disprel}
{A}_{3n,3n}^{} = 
\begin{pmatrix} 
 a_{1}^{} {A}_{1}^{t} (\frac{1}{\alpha},\frac{1}{\beta}) + a_{2}^{} {A}_{2}^{t} (\frac{1}{\alpha},\frac{1}{\beta}) + p {A}_{3}^{} (\alpha,\beta) 
& {A}_{1}^{}(\alpha,\beta)
& {A}_{2}^{}(\alpha,\beta) \\*[0.4cm]
 \nu {A}_{1}^{t} (\frac{1}{\alpha},\frac{1}{\beta})  
& I_{n,n}^{}
& O_{n,n}^{} \\*[0.4cm]
 \nu {A}_{2}^{t} (\frac{1}{\alpha},\frac{1}{\beta})
& O_{n,n}^{}
& I_{n,n}^{}
\end{pmatrix}, 
\end{equation}
where $O_{n,n}^{}$ and $I_{n,n}^{}$ are the $n \times n$ null and identity matrices, respectively.
After long and tedious algebra, the $n \times n$ matrices ${A}_{1}^{}, {A}_{2}^{}$ and ${A}_{3}^{}$ are computed and shown in~\cref{app-def-submat-fourier-NC,app-def-submat-fourier-DG} for the $P^{NC}_{1}$ and $P^{DG}_{1}$ methods, respectively.
\subsection{The dispersion relations}
\label{Disp-Rel}
In order to obtain non trivial solutions for the Fourier amplitudes $\widehat{\mathbf{X}}$ in~\eqref{disprel}, the so-called dispersion relation or characteristic polynomial $P_{n}^{} (\omega)$, namely, a relationship between the frequency $\omega$ and the wave numbers $k$ and $l$, is computed using the MAPLE software by setting
\begin{equation}
\label{P3n-Gen-disp-Rel}
P_{n}^{} (\omega) := \det \Big( {A}_{3n,3n}^{} - i \omega\,{B}_{3n,3n}^{} \Big) = 0. 
\end{equation}
Consequently, $P_{n}^{} (\omega)$ is a polynomial of degree $n$ in $\omega$, whose coefficients depend on the parameters $k$, $l$, $h$, $p$, $\nu$, $a_{1}^{}$ and $a_{2}^{}$. 

On the other hand, the governing equation~\eqref{cont-eq} is linear with constant coefficients and when periodic boundary conditions are employed, the solution may be examined by considering the behaviour of a single Fourier mode. We then seek a solution of~\eqref{cont-eq} of the form $u=\widehat{u}\,e^{i(kx+ly)}\,e^{i\omega t}$, where $\widehat{u}$ is the Fourier amplitude. Substitution into~\eqref{cont-eq} leads to the continuous or analytical (an) frequency
\begin{equation}
\label{freq-an}
\omega_{}^{an} = \underbrace{-(a_{1}^{} k + a_{2}^{} l)}_{\textstyle := \omega_{adv}^{an} \mathstrut} 
+ \underbrace{i \nu (k^2+l^2)}_{\textstyle := \omega_{dif}^{an} \mathstrut}.
\end{equation}
The purpose of the subsequent analyses is to compare the $n$ discrete frequencies $\omega$ solutions of~\eqref{P3n-Gen-disp-Rel} with the continuous counterpart obtained in~\eqref{freq-an}. The asymptotic expansion of the discrete frequencies solutions of $P_{n}^{} (\omega)=0$ in~\eqref{P3n-Gen-disp-Rel}, for inﬁnitesimal mesh spacing $h$, are computed by using the MAPLE software for the $P_{1}^{NC}$ and $P_{1}^{DG}$ schemes. Unfortunately, MAPLE is unable to compute such asymptotic expansions in the general case, and we need to consider separately the case $\mathbf{a}=\mathbf{0}$ (\cref{mainthNC1,mainthNC2}) and the case $\nu=0$ and $p=0$ (\cref{mainth-NC-adv}). The asymptotic expansions show that a superconvergence result is obtained for the discrete frequency of the $P_{1}^{NC}$ scheme compared to the $P_{1}^{DG}$ one, denoted by $\omega_{1}^{nc}$ and $\omega_{1}^{dg}$, respectively.
\begin{theorem}
\label{mainthNC1}
In the absence of the advection term (for $\mathbf{a}=\mathbf{0}$), and in the limit as mesh spacing $h \rightarrow 0$, we obtain for all $p \in \mathbb{R}$, the asymptotic results for the frequencies $\omega_{1}^{nc}$ and $\omega_{1}^{dg}$ solutions of~\eqref{P3n-Gen-disp-Rel} for the $P_{1}^{NC}$ and $P_{1}^{DG}$ schemes, respectively, 
\begin{align*}
& \omega_{1}^{nc} = \omega_{dif}^{an} + \nu f_{1}^{}(k,l,p)\, \boldsymbol{i h^{4}_{}} + O(h^6),  \\
& \omega_{1}^{dg} = \omega_{dif}^{an} + \nu g_{1}^{}(k,l,p)\, \boldsymbol{i h^{2}_{}} + O(h^4),
\end{align*}
where $f_{1}^{}$ and $g_{1}^{}$ are polynomials in $k, l$ and $p$.

The additional frequencies $\omega_{2,3}^{nc}$ and $\omega_{2,3,...,6}^{dg}$ behave as high frequency modes of the form $\frac{i \nu}{h_{}^{2}} ( p + p_{j}^{}) + O(1)$, where $p_{j}^{} \in \left \{2, \frac{14}{3}\right \}$ and $p_{j}^{} \in \left \{4, \frac{14}{3}, 6, 8, 12 \right \}$ for the $P_{1}^{NC}$ and $P_{1}^{DG}$ schemes, respectively. Such frequencies hence decay exponentially fast provided that $p + p_{j}^{} > 0$ and they have no impact in practical computations.

Further, these results demonstrate that the discrete frequencies solutions of~\eqref{P3n-Gen-disp-Rel} are stable at leading order when $p > -2$ and $p > -4$ for the $P^{NC}_{1}$ and $P_{1}^{DG}$ methods, respectively.
\end{theorem}
\begin{theorem}
\label{mainthNC2}
There exists a unique value of $p$ in~\cref{mainthNC1}, namely $p=\frac{1}{3}$ for the $P_{1}^{NC}$ scheme and $p \simeq -1.3822831$ for the $P_{1}^{DG}$ scheme, for which these numerical schemes are stable and in the limit as mesh spacing $h \rightarrow 0$ the following asymptotic results are obtained
\begin{align*}
& \omega_{1}^{nc}  = \omega_{dif}^{an} - \nu \Big( \frac{1}{90} (k_{}^{4}+l_{}^{4}) -0.005952 (k_{}^{3}l_{}^{}+k_{}^{}l_{}^{3}) +  0.001984 k_{}^{2}l_{}^{2} \Big) \boldsymbol{kl} \,i h^{4}_{} + O(h^6), \\
& \omega_{1}^{dg}  = \omega_{dif}^{an} - \nu \big( 0.03169 (k_{}^{2}+l_{}^{2}) - 0.02511 kl \big) \boldsymbol{kl} \,i h^{2}_{} + O(h^4).
\end{align*}
Therefore, the choice $\boldsymbol{l=0}$ (or $\boldsymbol{k=0}$ by exchanging $k$ and $l$) yields
\begin{align*}
& \omega_{1}^{nc} = \omega_{dif}^{an} -0.0005864 \,\nu k_{}^{8}  \boldsymbol{i h^{6}_{}} + O(h^8), \qquad \mbox{for } p=\frac{1}{3} \\
& \omega_{1}^{dg} = \omega_{dif}^{an} -0.0010047 \,\nu k_{}^{6}  \boldsymbol{i h^{4}_{}} + O(h^6), \qquad \mbox{for } p  \simeq -1.3822831.
\end{align*}
\end{theorem}

The values $p=\frac{1}{3}$ for the $P_{1}^{NC}$ scheme and $p \simeq -1.3822831$ for the $P_{1}^{DG}$ scheme are derived algebraically in the purpose of cancelling appropriate terms in the asymptotic results of~\cref{mainthNC1}. In particular, the value $p \simeq -1.3822831$ is obtained by solving 
\begin{equation}
\frac{1}{4} p^5 + \frac{89}{12} p^4 + \frac{242}{3} p^3 + \frac{1199}{3} p^2 + 880\,p + 640 = 0.
\label{Choice-p-P1DG}
\end{equation}
The four other roots are smaller than $-4$ in~\eqref{Choice-p-P1DG} and then lead to unstable results.

\begin{figure}[h!]
\begin{center}
\begin{tabular}{c|c}
\label{Fig_Disp-rel-min-max}
$\omega_1^{nc}/\omega_{}^{an}$  &   $\omega_1^{dg}/ \omega_{}^{an}$  \\*[0.1cm] 
\hline
\begin{minipage}{48mm}
\begin{center}
\includegraphics*[width=47mm]{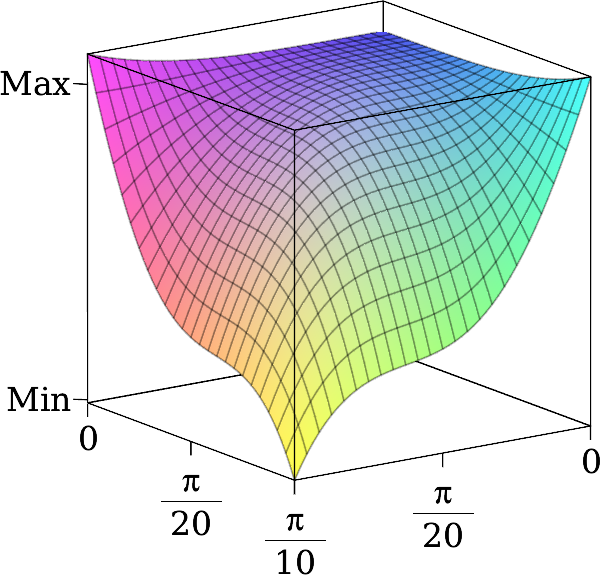}
\end{center}
\end{minipage} 
&
\begin{minipage}{48mm}
\begin{center}
\includegraphics*[width=47mm]{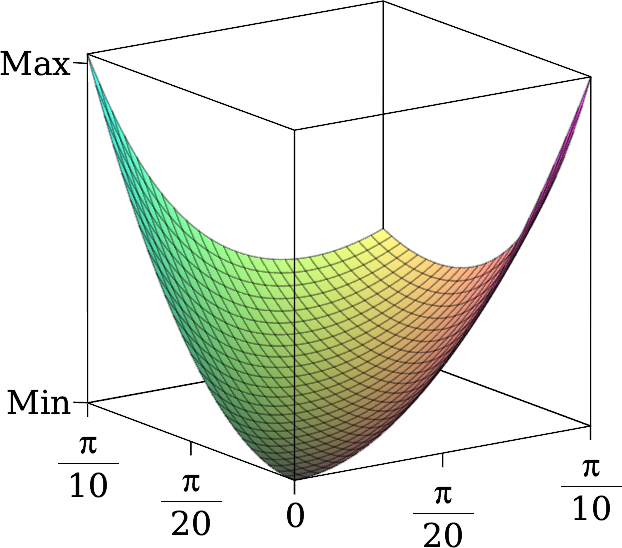}
\end{center}
\end{minipage}
\vrule height 25mm width 0 pt depth 23mm \\ 
$p=1$, \small{(Min,Max)=$(0.99996,1)$} & $p=1$, \small{(Min,Max)=$(1,1.0029)$} \\*[0.2cm]
\begin{minipage}{48mm}
\begin{center}
\includegraphics*[width=47mm]{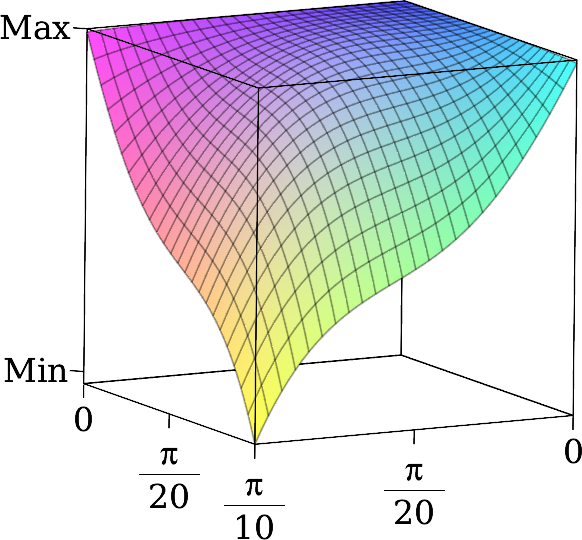}
\end{center}
\end{minipage} 
&
\begin{minipage}{48mm}
\begin{center}
\includegraphics*[width=48mm]{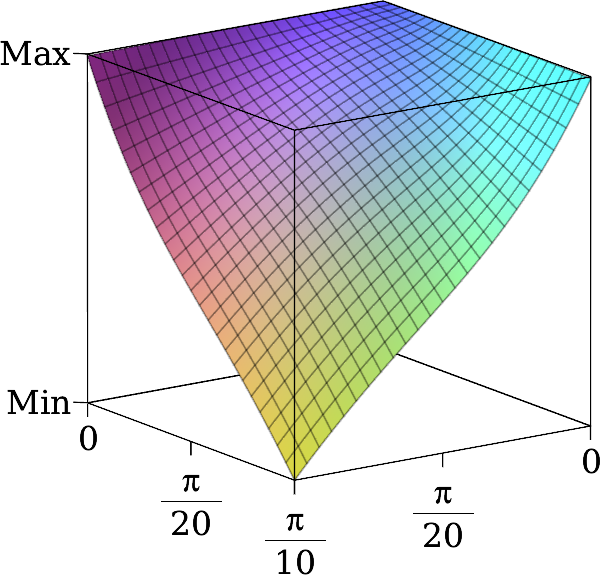}
\end{center}
\end{minipage}
\vrule height 25mm width 0 pt depth 23mm \\ 
$p=\frac{1}{3}$, \small{(Min,Max)=$(0.99995,1)$}  &  $p \simeq -1.382283$, \small{(Min,Max)=$(0.9981, 1)$} \\*[0.2cm]
\end{tabular}
\caption{The normalized frequencies $\omega_1^{nc}/\omega_{}^{an}$ and $\omega_1^{dg}/ \omega_{}^{an}$, for $0 \leq kh, lh \leq \pi/10$, in the case $p=1$ and $p=1/3$ for the $P_{1}^{NC}$ scheme, and $p=1$ and $p=-1.382283$ for the $P_{1}^{DG}$ method. The maximum and minimum (Max, Min) values of the normalized frequencies are specified in the figure.}
\end{center}
\end{figure}

The normalized frequencies $\omega_{1}^{nc}/\omega_{dif}^{an}$ and $\omega_{1}^{dg}/\omega_{dif}^{an}$ are graphed in~\cref{Fig_Disp-rel-min-max} in the case $0 \leq kh, lh \leq \pi/10$, and their maximum and minimum values are specified for specific values of $p$. Such results illustrate the very accurate behaviour of the $P^{NC}_{1}$ method observed in~\cref{mainthNC1} and~\cref{mainthNC2}, at least for small wavenumbers. 
\begin{theorem}
\label{mainth-NC-adv}
In the absence of the diffusion term, namely in the case $\nu=0$ and $p=0$, and in the limit as mesh spacing $h \rightarrow 0$, we obtain the asymptotic results for the frequencies $\omega_{1}^{nc}$ and $\omega_{1}^{dg}$ solutions of~\eqref{P3n-Gen-disp-Rel} for the $P_{1}^{NC}$ and $P_{1}^{DG}$ schemes, respectively,
\begin{align*}
& \omega_{1}^{nc} =  \omega_{adv}^{an} + \Big(\frac{1}{180} \left( a_{1}^{} k^{5}_{} + a_{2}^{} l^{5}_{} \right) - \frac{1}{72} \,kl \left(k-l \right) \left( a_{1}^{} k^{2}_{} - a_{2}^{} l^{2}_{} \right) \Big) \boldsymbol{h_{}^{4}} + \boldsymbol{i}\,O(\boldsymbol{h_{}^{5}}), \\
& \omega_{1}^{dg} = \omega_{adv}^{an} + \frac{1}{144 \left(a_{1}^{} \! + a_{2}^{} \right)} \left( k^{2}_{}\! -kl \! + l^{2}_{} \right) \Big( (a_{1}^{} k - a_{2}^{} l)^{2}_{} + 2 a_{1}^{} a_{2}^{} (k^{2}_{}\! + l^{2}_{}) \Big) \boldsymbol{i h_{}^{3}} \! + O(h_{}^{4}), 
\label{HFM-NC}
\end{align*}
The additional frequencies $\omega_{2,3}^{nc}$ and $\omega_{2,3,...,6}^{dg}$ behave as high frequency modes of the form $\frac{c_{j}^{} + i\,d_{j}^{}}{h} + O(1)$, where $c_{j}^{}>0$ and $d_{j}^{}=0$, $j=2,3,$ for the $P_{1}^{NC}$ scheme, and $\mid \!\!c_{j}^{} \!\!\mid \,\leq 9.863 \,\mbox{ and }\, 6.690 \leq d_{j}^{}  \leq 8.655,$ $j=2,3,\cdots,6$, for the $P_{1}^{DG}$ method. Such frequencies are hence stable, they decay exponentially fast and have no impact in practical computations. 
\end{theorem}
\section{Conclusion}
\label{se:con}
This paper appears to be the ﬁrst study of the Fourier/stability analysis of the 2D advection-diffusion equation for the linear discontinuous $P_{1}^{DG}$ and non conforming $P_{1}^{NC}$ Galerkin discretizations. The discrete formulations are stabilized by means of a family of numerical ﬂuxes and the dispersion relations are constructed. A long time stability result is proven for the $P_{1}^{NC}$ and $P_{1}^{DG}$ schemes and ﬂuxes examined here by inspecting the discrete frequencies, solutions of the dispersion relations. A super-convergent result is also demonstrated for the $P_{1}^{NC}$ discrete frequency ($O(h^4)$ accuracy) compared to the $P_{1}^{DG}$ one ($O(h^2)$ accuracy). Finally, these results are in line with those obtained in \cite{LeRoux-JCP-2024} for the shallow-water equations.
\appendix
\renewcommand{\thesection}{\Alph{section}}
\section{The vectors \boldmath$\Lambda^{}_{1}, \Lambda^{}_{2}, \Lambda^{}_{3}, \Lambda^{}_{4},$ for the \boldmath$P_{1}^{NC}$ scheme}
\label{app-def-submatNC}
{\small
\begin{equation*}
\label{Disc-oper-DG-2}
\boldsymbol{\Lambda}^{}_{1} = 
\begin{pmatrix}  
 2 u^{}_{1} -   u^{}_{3} + u^{}_{10} - u^{}_{11} - u^{}_{14} - u^{}_{15} + u^{}_{16}            \\
-u^{}_{1} + u^{}_{2} -   u^{}_{3} + u^{}_{4} - u^{}_{10} + u^{}_{11}                  \\
 u^{}_{2} +  u^{}_{3} - u^{}_{10} - u^{}_{11} + u^{}_{12} - u^{}_{13}                         
\end{pmatrix},
\end{equation*}
\begin{equation*}
\boldsymbol{\Lambda}^{}_{2} = 
\begin{pmatrix}  
u^{}_{1} + u^{}_{2} - u^{}_{14} - u^{}_{15} + u^{}_{16} - u^{}_{17} \\
-u^{}_{1} + u^{}_{2} - u^{}_{3} + u^{}_{5} + u^{}_{15} - u^{}_{17} \\
-u^{}_{1} + 2 u^{}_{3} - u^{}_{11} + u^{}_{12} - u^{}_{13}- u^{}_{15} + u^{}_{17} 
\end{pmatrix}.
\end{equation*}
\begin{equation*}
\boldsymbol{\Lambda}^{}_{3} = 
\begin{pmatrix}  
-4 u^{}_{1} + u^{}_{2} + u^{}_{3} + u^{}_{4} - u^{}_{5} - u^{}_{10} + u^{}_{11} 
+ u^{}_{14} + u^{}_{15} - u^{}_{16} + u^{}_{17} - u^{}_{18} + u^{}_{19} \\*[0.1cm]
u^{}_{1} - 4 u^{}_{2} + u^{}_{3} + u^{}_{4} + u^{}_{5} - u^{}_{6} + u^{}_{7} + u^{}_{10} 
- u^{}_{11} - u^{}_{15} + u^{}_{17} + u^{}_{20} - u^{}_{21} \\*[0.1cm]
u^{}_{1} + u^{}_{2} -4 u^{}_{3} - u^{}_{4} + u^{}_{5} + u^{}_{8} - u^{}_{9} 
+ u^{}_{10} + u^{}_{11} - u^{}_{12} + u^{}_{13} + u^{}_{15} - u^{}_{17}  
\end{pmatrix},
\end{equation*}
\begin{equation*}
\label{Disc-oper-DG-122}
\boldsymbol{\Lambda}^{}_{4} =
\begin{pmatrix}  
-2\theta^{x}_{1} + \theta^{x}_{2} + \theta^{x}_{4} - \theta^{x}_{5} + \theta^{x}_{17} - \theta^{x}_{18 } + \theta^{x}_{19} - \theta^{y}_{1} + \theta^{y}_{2} + \theta^{y}_{3} + \theta^{y}_{4} - \theta^{y}_{5} - \theta^{y}_{15} \\*[0.1cm]
-\theta^{x}_{2} - \theta^{x}_{3} + \theta^{x}_{4} + \theta^{x}_{5} + \theta^{x}_{20} - \theta^{x}_{21} - \theta^{y}_{1} - \theta^{y}_{2} + \theta^{y}_{4} + \theta^{y}_{5} - \theta^{y}_{6} + \theta^{y}_{7} \\*[0.1cm]
\theta^{x}_{1} + \theta^{x}_{2} - \theta^{x}_{3} - \theta^{x}_{4} + \theta^{x}_{5} - \theta^{x}_{11} + \theta^{y}_{2} - 2\theta^{y}_{3} - \theta^{y}_{4} + \theta^{y}_{5} + \theta^{y}_{8} - \theta^{y}_{9} + \theta^{y}_{10}
\end{pmatrix},
\end{equation*}
}
\section{The vectors \boldmath$\Lambda^{}_{1}, \Lambda^{}_{2}, \Lambda^{}_{3}, \Lambda^{}_{4},$ for the \boldmath$P_{1}^{DG}$ scheme}
\label{app-def-submatDG}
{\small
\begin{equation*}
\boldsymbol{\Lambda}^{}_{1} = 2
\begin{pmatrix}  
 u^{}_{1} -   u^{}_{3} +   u^{}_{11} - u^{}_{12}             \\
-u^{}_{1} + 2 u^{}_{2} -   u^{}_{3}  - u^{}_{11} + u^{}_{12} \\
 u^{}_{2} + 2 u^{}_{3} - 3 u^{}_{10}                         \\
-u^{}_{1} +   u^{}_{3} + 2 u^{}_{4}  - u^{}_{5}  - u^{}_{6}  \\
 u^{}_{1} -   u^{}_{3} +   u^{}_{5}  - u^{}_{6}              \\
-3 u^{}_{2} + u^{}_{4} + 2 u^{}_{6} 
\end{pmatrix}\!,
\,
\boldsymbol{\Lambda}^{}_{2} = 2
\begin{pmatrix}  
2 u^{}_{1} +   u^{}_{2} - 3 u^{}_{8}             \\
-u^{}_{1} + 2 u^{}_{2} -   u^{}_{3}  - u^{}_{7} + u^{}_{9} \\
-u^{}_{1} +  u^{}_{3}  + u^{}_{7} - u^{}_{9}                         \\
-u^{}_{1} +   u^{}_{3} + u^{}_{4} - u^{}_{6}  \\
 u^{}_{1} -   u^{}_{3} - u^{}_{4} + 2 u^{}_{5}  - u^{}_{6}              \\
-3 u^{}_{2} + u^{}_{5} + 2 u^{}_{6} 
\end{pmatrix}\!.
\end{equation*}
\begin{equation*}
\boldsymbol{\Lambda}^{}_{3} = 2
\begin{pmatrix}  
-5 u^{}_{1} + u^{}_{2} + u^{}_{3} + u^{}_{4} - u^{}_{5} + 3 u^{}_{8} - u^{}_{11} + u^{}_{12}  \\
u^{}_{1} - 5 u^{}_{2} + u^{}_{3} + 3 u^{}_{6} + u^{}_{7} - u^{}_{9} + u^{}_{11} - u^{}_{12}   \\
u^{}_{1} + u^{}_{2} -5 u^{}_{3} - u^{}_{4} + u^{}_{5} - u^{}_{7} + u^{}_{9} + 3 u^{}_{10}     \\
u^{}_{1} - u^{}_{3}-5 u^{}_{4} + u^{}_{5} + u^{}_{6} + 3 u^{}_{14} + u^{}_{16} - u^{}_{17}   \\
-u^{}_{1} + u^{}_{3} + u^{}_{4} -5 u^{}_{5} + u^{}_{6} + u^{}_{13} - u^{}_{15} + 3 u^{}_{18}  \\
3u^{}_{2} + u^{}_{4} + u^{}_{5} - 5 u^{}_{6} - u^{}_{13} + u^{}_{15} - u^{}_{16} + u^{}_{17} 
\end{pmatrix},
\end{equation*}
\begin{equation*}
\boldsymbol{\Lambda}^{}_{4} = 2
\begin{pmatrix}  
-\theta^{x}_{1} + \theta^{x}_{2} + \theta^{x}_{4} - \theta^{x}_{5} - 2 \theta^{y}_{1} + \theta^{y}_{2} + \theta^{y}_{3} + \theta^{y}_{4} - \theta^{y}_{5} \\
-2 \theta^{x}_{2} - \theta^{x}_{3} + 3 \theta^{x}_{6} - \theta^{y}_{1} - 2 \theta^{y}_{2} + 3 \theta^{y}_{6} \\
\theta^{x}_{1} + \theta^{x}_{2} - 2 \theta^{x}_{3} - \theta^{x}_{4} + \theta^{x}_{5} + \theta^{y}_{2} - \theta^{y}_{3} - \theta^{y}_{4} + \theta^{y}_{5} \\
-2 \theta^{x}_{4} - \theta^{x}_{6} + 3 \theta^{x}_{14} - \theta^{y}_{4} + \theta^{y}_{5} + \theta^{y}_{16} - \theta^{y}_{17} \\
\theta^{x}_{4} - \theta^{x}_{5} + \theta^{x}_{13} - \theta^{x}_{15} - 2 \theta^{y}_{5} - \theta^{y}_{6} + 3 \theta^{y}_{18} \\
\theta^{x}_{4} + \theta^{x}_{5} - 2 \theta^{x}_{6} - \theta^{x}_{13} + \theta^{x}_{15} + \theta^{y}_{4} + \theta^{y}_{5} - 2 \theta^{y}_{6} - \theta^{y}_{16} + \theta^{y}_{17} 
\end{pmatrix}.
\end{equation*}
}
\section{The matrices \boldmath${A}^{}_{1}, {A}^{}_{2}, {A}^{}_{3}$, for the \boldmath$P_{1}^{NC}$ scheme}
\label{app-def-submat-fourier-NC}
{\small
$$\begin{array}{llll}
\alpha^{}_{1} := \dfrac{\alpha}{\beta}+\dfrac{\beta}{\alpha}, & \quad \alpha^{}_{3} := \alpha+\dfrac{1}{\alpha}, & \quad \alpha^{}_{5} := \beta+\dfrac{1}{\beta}, & \quad \alpha^{}_{7} := -2-\dfrac{\alpha^{2}_{}}{\beta^{2}_{}}-\beta^{2}_{}, \\*[0.3cm]
\alpha^{}_{2} := \dfrac{\alpha}{\beta}-\dfrac{\beta}{\alpha}, & \quad \alpha^{}_{4} := \alpha-\dfrac{1}{\alpha}, & \quad \alpha^{}_{6} := \beta-\dfrac{1}{\beta}, & \quad \alpha^{}_{8} := -2-\dfrac{\beta^{2}_{}}{\alpha^{2}_{}}-\alpha^{2}_{},
\end{array}$$
\begin{displaymath}
\setlength\arraycolsep{2pt}
{A}_{1}^{}(\alpha,\beta) = \dfrac{1}{h}
{\small{
\begin{pmatrix}
\alpha^{}_{7}   & \alpha \alpha^{}_{1}  & \alpha \alpha^{}_{5} \\*[0.1cm]
\alpha \alpha^{}_{1} & -1-\alpha^{2}_{} & \alpha^{}_{4}   \\*[0.1cm]
\alpha \alpha^{}_{5} & \alpha^{}_{4}    & -1-\alpha^{2}_{}
\end{pmatrix}, \quad
}}
{A}_{2}^{} (\alpha,\beta) =  \dfrac{1}{h}
{\small{
\begin{pmatrix}
-1-\beta^{2}_{} & \alpha^{}_{6}    & \beta \alpha^{}_{3} \\*[0.1cm]
\alpha^{}_{6}    & -1-\beta^{2}_{} & \beta \alpha^{}_{1} \\*[0.1cm]
\beta \alpha^{}_{3}  & \beta \alpha^{}_{1}  & \alpha^{}_{8}
\end{pmatrix}, \quad
}}
\end{displaymath}
\begin{displaymath}
\setlength\arraycolsep{2pt}
{A}_{3}^{}(\alpha,\beta)  =  \dfrac{1}{h_{}^{2}}
{\small{
\begin{pmatrix}
-\alpha^{2}_{1}-\alpha^{2}_{5} & \alpha^{}_{1} \alpha^{}_{3}    & \alpha^{}_{3}\alpha^{}_{5}     \\*[0.2cm]
\alpha^{}_{1} \alpha^{}_{3}            & -\alpha^{2}_{3}-\alpha^{2}_{5} & \alpha^{}_{1} \alpha^{}_{5}            \\*[0.2cm]
\alpha^{}_{3}\alpha^{}_{5}     & \alpha^{}_{1} \alpha^{}_{5}    & -\alpha^{2}_{1}-\alpha^{2}_{3}
\end{pmatrix}.
}}
\end{displaymath}
}
\section{The matrices \boldmath${A}^{}_{1}, {A}^{}_{2}, {A}^{}_{3}$, for the \boldmath$P_{1}^{DG}$ scheme}
\label{app-def-submat-fourier-DG}
{\small
$$\begin{array}{lll}
\beta^{}_{1} = 3- \dfrac{\beta^{2}_{}}{\alpha^{2}_{}}-\beta^{2}_{},    & \quad \beta^{}_{3} = 3-\dfrac{\alpha^{2}_{}}{\beta^{2}_{}}-\alpha^{2}_{}, 
& \quad \beta^{}_{5} = 3-\dfrac{\alpha^{2}_{}}{\beta^{2}_{}}-\dfrac{1}{\beta^{2}_{}}, \\*[0.3cm]
\beta^{}_{2} = 3-\dfrac{1}{\alpha^{2}_{}}-\dfrac{1}{\beta^{2}_{}}, & \quad \beta^{}_{4} = 3-\dfrac{\beta^{2}_{}}{\alpha^{2}_{}}-\dfrac{1}{\alpha^{2}_{}}, 
& \quad \beta^{}_{6} = 3-\alpha^{2}_{} - \beta^{2}_{}, 
\end{array}$$
\begin{displaymath}
\setlength\arraycolsep{2pt}
{A}_{1}^{}(\alpha,\beta) = \dfrac{2}{h}
{\small{
\begin{pmatrix}
-1 & \beta & 0 & \alpha \beta & -\beta^{2}_{} & 0 \\*[0.1cm] 
0 & -2 & -\dfrac{1}{\alpha} & 0 & 0 & 3 \\*[0.1cm] 
\dfrac{\alpha}{\beta} & \alpha & -2 & -\alpha^{2}_{} & \alpha \beta & 0 \\*[0.1cm] 
0 & 0 & 3 & -2 & 0 & -\dfrac{1}{\alpha} \\*[0.1cm]
-\dfrac{\alpha^{2}_{}}{\beta^{2}_{}} & \dfrac{\alpha^{2}_{}}{\beta} & 0 & \dfrac{\alpha}{\beta} & -1 & 0 \\*[0.1cm]
\dfrac{\alpha^{2}_{}}{\beta} & -\alpha^{2}_{} & 0 & \alpha & \beta & -2
\end{pmatrix}\!, 
{A}_{3}^{}(\alpha,\beta) = \dfrac{2}{h_{}^{2}}
{\small{
\begin{pmatrix}
-5 & \beta & \dfrac{\beta}{\alpha} & \alpha \beta &  \beta^{}_{1} & \dfrac{\beta}{\alpha^{2}_{}} \\*[0.2cm] 
\dfrac{1}{\beta} & -5 & \dfrac{1}{\alpha} & \dfrac{\alpha}{\beta^{2}_{}} & \dfrac{\beta}{\alpha^{2}_{}}  &  \beta^{}_{2} \\*[0.2cm] 
\dfrac{\alpha}{\beta} & \alpha & -5 & \beta^{}_{3} &  \alpha \beta &  \dfrac{\alpha}{\beta^{2}_{}} \\*[0.2cm] 
\dfrac{1}{\alpha \beta} & \dfrac{\beta^{2}_{}}{\alpha} & \beta^{}_{4} & -5 &  \dfrac{\beta}{\alpha} & \dfrac{1}{\alpha}   \\*[0.2cm]
\beta^{}_{5} & \dfrac{\alpha^{2}_{}}{\beta} & \dfrac{1}{\alpha \beta} & \dfrac{\alpha}{\beta} &  -5 &  \dfrac{1}{\beta} \\*[0.2cm]
\dfrac{\alpha^{2}_{}}{\beta} & \beta^{}_{6} & \dfrac{\beta^{2}_{}}{\alpha} & \alpha &  \beta &  -5
\end{pmatrix}\!\!, 
}}
}}
\end{displaymath}
\begin{displaymath}
\setlength\arraycolsep{2pt}
{A}_{2}^{}(\alpha,\beta) = \dfrac{2}{h}
{\small{
\begin{pmatrix}
-2 & \beta & \dfrac{\beta}{\alpha} & \alpha \beta & -\beta^{2}_{} & 0 \\*[0.1cm] 
-\dfrac{1}{\beta} & -2 & 0 & 0 & 0 & 3 \\*[0.1cm]
0 & \alpha & -1 & -\alpha^{2}_{} & \alpha \beta & 0 \\*[0.1cm]
0 & \dfrac{\beta^{2}_{}}{\alpha} & -\dfrac{\beta^{2}_{}}{\alpha^{2}_{}} & -1 & \dfrac{\beta}{\alpha} & 0 \\*[0.1cm]
3 & 0 & 0 & 0 & -2 & -\dfrac{1}{\beta} \\*[0.1cm]
0 & -\beta^{2}_{} & \dfrac{\beta^{2}_{}}{\alpha} & \alpha & \beta & -2 
\end{pmatrix}.
}}
\end{displaymath}
}
%
%
%
%
%

\end{document}